\documentstyle[amssymb,10pt]{amsart}%\input amssym.def
\newtheorem{theorem}{Theorem}[section]
\newtheorem{lemma}[theorem]{Lemma}
\newtheorem{proposition}[theorem]{Proposition}
 
\newtheorem{definition}[theorem]{Definition}

\theoremstyle{remark}

\newcommand{\op}{{\text{op}}}
 
\newcommand{\nc}{\newcommand}
\nc{\Symm}{{\on{Sym}}}
\nc{\Perm}{{\on{Perm}}}
\nc{\Iff}{\Leftrightarrow}
\newcommand{\on}{\operatorname}   
\newcommand{\eps}{\varepsilon}

 \nc{\cE}{{\cal E}}
\nc{\cL}{{\cal L}}

\nc{\SL}{{\mathfrak sl}}
\nc{\gt}{{\mathfrak gt}}
\nc{\grt}{{\mathfrak grt}}
\nc{\gtm}{{\mathfrak gtm}}
\nc{\grtm}{{\mathfrak grtm}}
\nc{\gtmd}{{\mathfrak gtmd}}
\nc{\grtmd}{{\mathfrak grtmd}}

\renewcommand{\a}{{\mathfrak{a}}}

\nc{\HH}{{\mathfrak H}}

\newcommand{\ul}{\underline}

\nc{\wh}{\widehat}\nc{\wt}{\widetilde}

\newcommand{\kk}{{\bf k}}

\newcommand{\ben}{\begin{enumerate}}
\newcommand{\een}{\end{enumerate}}

\newcommand{\cO}{{\mathcal O}}

\newcommand{\cS}{{\mathcal S}}

\begin{document}

\title{Quantization of $\Gamma$-Lie bialgebras}

\begin{abstract}
We introduce the notion of $\Gamma$-Lie bialgebras, where $\Gamma$
is a group. These objects give rise to cocommutative 
co-Poisson bialgebras, for which we construct quantization functors.
This enlarges the class of co-Poisson algebras for which a quantization 
is known. Our result relies on our earlier work, where we showed that 
twists of Lie bialgebras can be quantized; we complement this work 
by studying the behavior of this quantization under compositions of 
twists. 
\end{abstract}

\author{Benjamin Enriquez}
\address{IRMA (CNRS), rue Ren\'e Descartes, F-67084 Strasbourg, France}
\email{enriquez@@math.u-strasbg.fr}

\author{Gilles Halbout}
\address{IRMA (CNRS), rue Ren\'e Descartes, F-67084 Strasbourg, France}
\email{halbout@@math.u-strasbg.fr}

\maketitle 

We work over a field $\kk$ of characteristic $0$. 

\section{Introduction}

Recall that a co-Poisson bialgebra is a quadruple $(U,m,\Delta_0,\delta)$, 
where $(U,m,\Delta_0)$ is a cocommutative bialgebra and  
$\delta: U \to \wedge^2(U)$ is a derivation (for $m$), 
a coderivation (for $\Delta_0$) and satisfies the co-Jacobi identity.
A quantization of $(U,m,\Delta_0,\delta)$ is a bialgebra 
$(U_\hbar,m_\hbar,\Delta_\hbar)$ such that $U_\hbar \simeq
U[[\hbar]]$, 
$m_\hbar=m + O(\hbar)$, $\Delta_\hbar=\Delta + O(\hbar)$,
$\Delta_\hbar(a)-\Delta_\hbar^\op(a)=\hbar \delta(a)+ O(\hbar^2)$
($\Delta_\hbar^\op$ is the opposite coproduct).

If $(\a,\mu_\a)$ is a Lie algebra ($\mu_\a : 
\wedge^2(\a) \to \a$ is the Lie bracket), the co-Poisson 
bialgebra structures on $U(\a)$ correspond bijectively to the maps
$\delta_\a: \a\to \wedge^2(\a)$ such that $(\a,\mu_\a,\delta_\a)$
is a Lie bialgebra. The quantization of these co-Poisson 
bialgebras was obtained in \cite{EK}. 

To a triple $(\Gamma,\a,\theta_a)$, where $\Gamma$ is a group, 
$\a$ is a Lie algebra
and $\theta_\a:\Gamma \to \on{Aut}(\a,\mu_\a)$ is an action of 
$\Gamma$ on $\a$, one associates the $\Gamma$-graded cocommutative 
bialgebra $U(\a)\rtimes \Gamma$. The $\Gamma$-graded co-Poisson 
bialgebra structures on $U(\a)\rtimes\Gamma$ correspond bijectively 
to pairs $(\delta_\a,f)$, where $\delta_\a:\a\to \wedge^2(\a)$
is such that $(\a,\mu_\a,\delta_\a)$ is a Lie bialgebra, and 
$f:\Gamma\to \wedge^2(\a)$ satisfies some conditions (see Section 
\ref{gammaLie}); 
in particular, $f(\gamma)$ is a twist of $(\a,\mu_\a,\delta_\a)$ 
for any $\gamma\in\Gamma$. We call the resulting 5-uple
$(\a,\mu_\a,\delta_\a,\theta_\a,f)$ a $\Gamma$-Lie bialgebra.  
The main result of this paper is the quantization of the 
corresponding co-Poisson bialgebra structures. 

Examples of $\Gamma$-Lie bialgebras arise from the following situation: 
$G$ is a Poisson-Lie group with Lie bialgebra $(\a,\mu_\a,\delta_\a)$, 
and $\Gamma \subset G$ is a discrete subgroup. Another example is when $\a$ 
is a Kac-Moody Lie algebra $\a$, and $\Gamma$ is the extended Weyl group 
of $\a$. In the latter case, a quantization in known (\cite{MS}). 

\smallskip

To achieve our goal, we complement a result obtained in \cite{EH}, namely  
the compatibility of Etingof-Kazhdan (EK) quantization functors with 
twists of Lie bialgebras; this result is based on an alternative 
construction of these quantization functors (\cite{En}). We 
describe the behavior of this quantization under
composition of twists (Section \ref{comp:twists}). 

To give an idea
of the result of \cite{EH}, we formulate its main consequence: 
let $Q : \{$Lie bialgebras$\} \to 
\{$quantized universal enveloping (QUE) algebras$\}$, 
$(\a,\mu_\a,\delta_\a) = \a \mapsto Q(\a) = (Q(\a),m(\a),
\Delta(\a),1_\a,\eps_\a)$ be a 
quantization 
functor; to each classical twist $f_\a$ of $\a$ (i.e., 
$f_\a\in\wedge^2(\a)$ and $(\delta_\a \otimes
\on{id}_\a)(f_\a) + [f_\a^{13},f_\a^{23}] + $ cyclic permutations = 0), 
one associates $\on{F}(\a,f_\a)\in Q(\a)^{\otimes 2}$, with the 
following properties: (a) it is a cocycle for $Q(\a)$, i.e.,\footnote{We 
denote by $*$ the product in $Q(\a)^{\otimes k}$ for $k\geq 1$} 
$(\on{F}(\a,f_\a)\otimes 1_{\a}) * (\Delta(\a) \otimes 
\on{id})(\on{F}(\a,f_\a)) = (1_{\a} \otimes \on{F}(\a,f_\a)) * (\on{id} 
\otimes \Delta(\a))(\on{F}(\a,f_\a))$,
$(\varepsilon_{\a} \otimes \on{id})(\on{F}(\a,f_\a))  
= (\on{id} \otimes \varepsilon_{\a})(\on{F}(\a,f_\a)) = 1_{\a}$; 
this implies that\footnote{If $A$ is an algebra and $u\in A$
is invertible, $\on{Ad}(u) : A \to A$ is $x\mapsto uxu^{-1}$} 
${}^{\on{F}(\a,\delta_\a)}Q(\a):= (Q(\a),m(\a),\on{Ad}(\on{F}(\a,f_\a)) 
\circ \Delta(\a),1_{\a},\varepsilon_{\a})$ is a QUE algebra; 
(b) we have an isomorphism $\on{i}(\a,f_\a) : 
{}^{\on{F}(\a,f_\a)}Q(\a) \to Q(\a_{f_\a})$ of QUE algebras
(here $\a_{f_{\a}} = (\a,\mu_\a,\delta_\a + \on{ad}(f_\a))$, where $\on{ad}
(f_\a) : \a\to \wedge^2(\a)$ is $x\mapsto [f_\a,x\otimes 1 + 1 \otimes x]$).  

In Section \ref{comp:twists}, we study the behavior of 
the assignment $(\a,f_\a) \mapsto 
\on{F}(\a,f_\a)$ under the composition of twists. A composition of 
twists is a pair $(f_\a,f'_\a)$ such that $f_\a$ is a twist of $\a$, 
and $f'_\a$ is a twist of $\a_{f_\a}$. We formulate the main 
consequence of our results: 
(a) there exists an invertible $\on{v}(\a,f_\a,f'_\a)
\in Q(\a)$, such that $\on{F}(f_\a + f'_\a) = 
\on{v}(\a,f_\a,f'_\a)^{\otimes 2} * \on{i}(\a,f_\a)^{-1}
(\on{F}(\a_{f_{\a}},f'_\a)) * F(\a,f_\a) * 
\Delta(\a)(\on{v}(\a,f_\a,f'_\a))^{-1}$
(Theorem \ref{thm:i:v}); and (b) if $(f_\a,f'_\a,f''_\a)$ are 
such that $f'_\a$ is a twist of $\a_{f_\a}$ and $f''_\a$ is a twist 
of $\a_{f_\a}$, then $f_\a + f'_\a$ is a twist of $\a$, and 
$\on{v}(\a,f_\a+f'_\a,f''_\a) * \on{v}(\a,f_\a,f'_\a) 
= \on{v}(\a,f_\a,f'_\a+f''_\a) *
\on{i}(\a,f_\a)^{-1}(\on{v}(\a_{f_\a},f'_\a,f''_\a))$
(Theorem \ref{thm:i:v}).

\smallskip

We use these results in Section \ref{quantgamma} to construct a quantization 
of the $\Gamma$-graded co-Poisson bialgebras $U(\a)\rtimes \Gamma$ 
as $\Gamma$-graded bialgebras. This quantization is based on 
the facts that that for $\gamma\in\Gamma$, $f(\gamma)$ is a twist of $\a$, 
and for any $\gamma,\gamma'\in\Gamma$, 
$(f(\gamma),\wedge^2(\theta_\a(\gamma))(f(\gamma')))$ is a composition of 
twists for $\a$; we then use the results of \cite{EH} on quantization of 
twists and those of Section \ref{comp:twists} on their composition.

\smallskip

This paper is organized as follows. 
In Section \ref{gammaLie}, we define $\Gamma$-Lie bialgebras, 
the corresponding co-Poisson cocommutative bialgebras, and the 
problem of their quantization. 
In Section \ref{back}, we recall the formalism of 
(quasi-multi-bi)props, which is the natural framework of 
the approach of \cite{En} to quantization functors and of 
the results of \cite{EH} on quantization of twists.
In Section \ref{comp:twists}, we describe the behavior of 
composition of twists under quantization (Theorem \ref{thm:i:v} and 
Theorem \ref{thm:v:v}).
In Section \ref{quantgamma}, we apply these results to
the construction of quantizations of $\Gamma$-Lie bialgebras.

\section{$\Gamma$-Lie bialgebras}\label{gammaLie}

\subsection{$\Gamma$-Lie algebras and equivalent categories}

Define a group Lie algebra as a triple $(\Gamma,\a,\theta_\a)$, 
where $\Gamma$ is a group, $\a$ is a Lie algebra and 
$\theta_\a : \Gamma \to \on{Aut}(\a)$ is a group morphism. 
Group Lie algebras form a category, where a morphism 
$(\Gamma,\a,\theta_\a) \to (\Gamma',\a',\theta_{\a'})$ is the data 
of a group morphism $i_{\Gamma\Gamma'} : \Gamma \to \Gamma'$
and a Lie algebra morphism $i_{\a\a'} : \a\to \a'$, such that 
$i_{\a\a'}(\theta_{\a,\gamma}(x)) = \theta_{\a,i_{\Gamma\Gamma'}(\gamma)}
(i_{\a\a'}(x))$. 

If $\Gamma$ is a group, a $\Gamma$-Lie algebra is a pair
$(\a,\theta_\a)$, such that $(\Gamma,\a,\theta_\a)$ is a group 
Lie algebra. $\Gamma$-Lie algebras form a subcategory of group 
Lie algebras, where the morphisms are restricted by the condition 
$i_{\Gamma\Gamma} = \on{id}_{\Gamma}$. 

Define a group cocommutative bialgebra as a triple $(\Gamma,A,i)$, where 
$\Gamma$ is a group, $A$ is a cocommutative bialgebra, 
$A = \oplus_{\gamma\in\Gamma} A_\gamma$ is a decomposition of $A$, and 
$i : \kk\Gamma \to A$ is a bialgebra morphism, such that 
$A_\gamma A_{\gamma'}\subset A_{\gamma\gamma'}$, $\Delta_A(A_{\gamma}) 
\subset A_\gamma^{\otimes 2}$, and $i$ is compatible with the $\Gamma$-grading. 
A morphism $(\Gamma,A,i)\to (\Gamma',A',i')$ is the data of a group morphism 
$i_{\Gamma\Gamma'} : \Gamma\to \Gamma'$ and a bialgebra morphism 
$i_{AA'} : A\to A'$, such that $i_{AA'}(A_\gamma) \subset 
A'_{i_{\Gamma\Gamma'}(\gamma')}$, and $i_{AA'} \circ i = i' \circ
i_{\kk\Gamma,\kk\Gamma'}$ (where $i_{\kk\Gamma,\kk\Gamma'} : \kk\Gamma 
\to \kk\Gamma'$ is the morphism induced by $i_{\Gamma\Gamma'}$). 

We then define a $\Gamma$-cocommutative bialgebra as a pair $(A,i)$, such that 
$(\Gamma,A,i)$ is a group cocommutative bialgebra. $\Gamma$-cocommutative 
bialgebras form a 
category, where as before $i_{\Gamma\Gamma} = \on{id}_{\Gamma}$. 

The category of group (resp., $\Gamma$-) cocommutative bialgebras contains 
as a full subcategory the category of group (resp., $\Gamma$-) 
universal enveloping algebras, where $(A,\Gamma,i)$ satisfies the
additional requirement that $A_e$ is a universal enveloping algebra. 

Define a group commutative bialgebra (in a symmetric monoidal category $\cS$)
as a triple $(\Gamma,\cO,j)$, where $\Gamma$ is a group, $\cO$ is a 
commutative algebra (in $\cS$) with a decomposition 
$\cO = \oplus_{\gamma\in\Gamma} \cO_\gamma$, 
such that $\cO_{\gamma}\cO_{\gamma'} = 0$ for $\gamma\neq\gamma'$, 
algebra morphisms $\Delta_{\gamma'\gamma''} : \cO_{\gamma'\gamma''} \to 
\cO_{\gamma'} \otimes \cO_{\gamma''}$, 
$\eta : \kk \to \cO_e$ and $\varepsilon :\cO_e \to \kk$, satisfying 
axioms such that when $\Gamma$ is finite, these morphisms add up to a 
bialgebra structure on $\cO$; and $j : \cO \to \kk^\Gamma$ is a morphism 
of commutative algebras, compatible with the $\Gamma$-gradings and the 
maps $\Delta_{\gamma'\gamma''}$ on both sides. We define 
$\Gamma$-commutative bialgebras as above. 

We define the category of group (resp., $\Gamma$-) formal series Hopf 
(FSH) algebras as a full subcategory of the category of 
group (resp., $\Gamma$-) commutative bialgebras in $\cS = \{$pro-vector 
spaces$\}$
by the condition the $\cO_e$ (or equivalently, each $\cO_\gamma$) is a formal 
series algebra. 

\begin{proposition}
1) We have (anti)equivalences of categories $\{$group Lie algebras$\} 
\leftrightarrow \{$group universal enveloping algebras$\}
\leftrightarrow\{$group FHS algebras$\}$ (the last map is an antiequivalence). 

2) If $\Gamma$ is a group, these (anti)equivalences restrict to 
$\{\Gamma$-Lie algebras$\} 
\leftrightarrow \{\Gamma$-universal enveloping algebras$\}
\leftrightarrow\{\Gamma$-FHS algebras$\}$. 
\end{proposition}

{\em Proof.} We denote the $\Gamma$-universal enveloping algebra 
corresponding to a $\Gamma$-Lie algebra $(\Gamma,\a,\theta_\a)$ as
$U(\a) \rtimes \Gamma$. It is isomorphic to $U(\a) \otimes \kk\Gamma$ 
as a vector space; 
if we denote by $x\mapsto [x]$, $\gamma\mapsto [\gamma]$ the natural maps 
$\a\to U(\a) \rtimes \Gamma$, $\Gamma \to U(\a) \rtimes \Gamma$, then  
the bialgebra structure of $U(\a) \rtimes \Gamma$ is given by 
$[\gamma][x][\gamma^{-1}] = [\theta_\gamma(x)]$, $[\gamma][\gamma'] =
[\gamma\gamma']$, $[e] = 1$, $[x][x'] - [x'][x] = [[x,x']]$, 
$\Delta([x]) = [x]\otimes 1 + 1 \otimes [x]$, 
$\Delta([\gamma]) = [\gamma] \otimes [\gamma]$.

When $\Gamma$ is finite, the corresponding $\Gamma$-FSH algebra is then 
$(U(\a) \rtimes \kk\Gamma)^*$, and in general, this is 
$\oplus_{\gamma\in\Gamma} (U(\a) \otimes \kk\gamma)^*$. 
One checks that these are (anti)equivalences of categories. 
For example, if $A$ is a group universal enveloping algebra, then one recovers
$\Gamma$ as $\{$group-like elements of $A\}$ and $\a$ as $\{$primitive 
elements of $A\}$.  \hfill \qed \medskip 

\subsection{$\Gamma$-Lie bialgebras and equivalent categories}

A group Lie bialgebra 
is a 5-uple $(\Gamma,\a, \theta_\a, \delta_\a,f)$ 
where $(\Gamma,\a,\theta_\a)$ is a group Lie algebra, 
$\delta_\a : \a\to\wedge^2(\a)$ is\footnote{We view 
$\wedge^2(V)$ as a  subspace of $V^{\otimes 2}$.} such that $(\a,\delta_\a)$
is a Lie bialgebra, and $f : \Gamma \to \wedge^2(\a)$
is a map $\gamma\mapsto f_\gamma$, such that: 
(a) $\wedge^2(\theta_\gamma) \circ \delta \circ 
\theta_\gamma^{-1}(x) = \delta(x) + [f_\gamma,x\otimes 1 + 1 \otimes x]$
for any $x\in \a$, (b) $f_{\gamma\gamma'} = f_\gamma +
\wedge^2(\theta_\gamma)(f_{\gamma'})$, and (c) $(\delta\otimes 
\on{id})(f_\gamma)
+ [f_\gamma^{1,3},f_\gamma^{2,3}]$ + cyclic permutations  $= 0$. 

Group Lie bialgebras form a category, where a morphism 
$(\Gamma,\a, \theta_\a, \delta_\a,f)\to 
(\Gamma',\a', \theta_{\a'}, \delta_{\a}',f')$ is a group Lie 
algebra morphism $(\Gamma,\a, \theta_\a)\to 
(\Gamma',\a', \theta_{\a'})$, such that $i_{\a\a'} : \a\to \a'$
is a Lie bialgebra morphism and $\wedge^2(i_{\a\a'})(f_\gamma) = 
f'_{i_{\Gamma\Gamma'}(\gamma)}$. When $\Gamma$ is fixed, one defines
the category of $\Gamma$-Lie bialgebras as above. 

A co-Poisson structure on a group cocommutative bialgebra 
$(\Gamma,A,i)$ is a co-Poisson structure $\delta_A : A \to \wedge^2(A)$, 
such that $\delta_A(A_\gamma) \subset \wedge^2(A_\gamma)$. 
Co-Poisson group cocommutative bialgebras form a category, where 
a morphism $(\Gamma,A,i,\delta_A) \to (\Gamma',A',i',\delta_{A'})$
is a morphism $(\Gamma,A,i) \to (\Gamma',A',i')$ of group 
cocommutative bialgebras, compatible with the co-Poisson structures. 
Co-Poisson group universal enveloping algebras form a full subcategory of 
the latter category. One defines the full subcategories of co-Poisson 
$\Gamma$-cocommutative bialgebras and co-Poisson $\Gamma$-enveloping 
algebras as above. 

A Poisson structure on a group commutative bialgebra 
$(\Gamma,\cO,j)$ is a Poisson bialgebra structure $\{-,-\}: 
\wedge^2(\cO) \to \cO$, 
such that $\{\cO_\gamma,\cO_\gamma\} \subset \cO_{\gamma}$ and 
$\{\cO_\gamma,\cO_{\gamma'}\} = 0$ if $\gamma\neq \gamma'$. 
Poisson group bialgebras form a category, and Poisson group FSH
algebras form a full subcategory when $\cS = \{$pro-vector spaces$\}$. 
One defines the full subcategories of Poisson $\Gamma$-bialgebras and 
Poisson $\Gamma$-FSH algebras as above.

\medskip 

{\bf Example.} Let $G$ be a Poisson-Lie (e.g., algebraic) group, let 
$\Gamma \subset G$ be a subgroup (which we view as an abstract group). 
We define $\theta_\gamma := \on{Ad}(\gamma)$, where $\on{Ad} : G \to 
\on{Aut}_{Lie}(\a)$ is the adjoint action. If $P : G \to \wedge^2(\a)$ is the 
Poisson bivector, satisfying $P(gg') = P(g') + \wedge^2(\on{Ad}(g))(P(g'))$, 
then we set $f_\gamma := -P(\gamma)$. Then $(\a,\Gamma,f)$ is a 
$\Gamma$-Lie bialgebra. 

\medskip 

{\bf Example.} Assume that $(\a,r_\a)$ is a quasitriangular Lie 
bialgebra and $\theta : \Gamma \to \on{Aut}(\a,t_\a)$ is an action of 
$\Gamma$ on $\a$ by Lie algebra automorphisms preserving $t_\a := r_\a +
r_\a^{2,1}$. If
we set $f_\gamma := \theta_\gamma^{\otimes 2}(r) - r$, then $(\a,\theta,f)$
is a $\Gamma$-Lie bialgebra (we call this a quasitriangular $\Gamma$-Lie
bialgebra). For example, $\a$ is a Kac-Moody Lie algebra, 
and $\Gamma = \tilde W$ is the extended Weyl group of $\a$.

\medskip

\begin{proposition}
1) We have category (anti)equivalences $\{$group bialgebras$\}
\leftrightarrow \{$co-Poisson group universal enveloping 
algebras$\}\leftrightarrow \{$Poisson group FSH algebras$\}$. 

2) These restrict to category (anti)equivalences
$\{\Gamma$-bialgebras$\}
\leftrightarrow \{$co-Poisson $\Gamma$-universal enveloping 
algebras$\}\leftrightarrow \{$Poisson $\Gamma$-FSH algebras$\}$. 
\end{proposition}

{\em Proof.} If $(\a,\theta_\a,\delta_\a)$ is a $\Gamma$-Lie bialgebra, 
then the co-Poisson structure on $A:= U(\a) \rtimes \Gamma$ is given by 
$\delta_A([x]) = [\delta_\a(x)]$, and $\delta_A([\gamma])
= -[f_\gamma]([\gamma]\otimes [\gamma])$. (Here we also denote by 
$x\mapsto [x]$ the natural map $\wedge^2(\a) \to \wedge^2(U(\a) \rtimes 
\Gamma$).) One checks that this establishes the desired (anti)equivalences. 
\hfill \qed \medskip

\subsection{The problem of quantization of $\Gamma$-Lie bialgebras}

Define a $\Gamma$-graded bialgebra (in a symmetric monoidal category $\cS$) 
as a bialgebra $A$ (in $\cS$), equipped with a grading 
$A = \oplus_{\gamma\in \Gamma} A_\gamma$, such that $A_\gamma A_{\gamma'}
\subset A_{\gamma\gamma'}$ and $\Delta_A(A_{\gamma}) \subset 
A_\gamma^{\otimes 2}$. 

Assume that $A$ is a $\Gamma$-graded bialgebra in the category of 
topologically free $\kk[[\hbar]]$-modules, quasicocommutative (in the sense 
that $A_0 := A/\hbar A$ is cocommutative). 
Then we get a co-Poisson structure on $A_0$. It is $\Gamma$-graded, in 
the sense that $\delta_{A_0}((A_0)_\gamma) \subset \wedge^2((A_0)_\gamma)$.
We therefore get a classical limit functor $\on{class} : 
\{\Gamma$-graded quasicocommutative 
bialgebras$\} \to \{\Gamma$-graded co-Poisson bialgebras$\}$. 

\begin{definition} A quantization functor for $\Gamma$-Lie bialgebras
is a functor $\{$co-Poisson $\Gamma$-universal enveloping 
algebras$\} \to \{\Gamma$-graded quasicocommutative bialgebras$\}$, 
right inverse to $\on{class}$. 
\end{definition}

We define the category of group-graded bialgebras as follows: objects are 
pairs $(\Gamma,A)$, where $\Gamma$ is a group and $A$ is a $\Gamma$-graded
bialgebra. A morphism $(\Gamma,A)\to (\Gamma',A')$ is the pairs of a group 
morphism $i_{\Gamma\Gamma'} : \Gamma \to \Gamma'$ and a bialgebra morphism 
$i_{AA'} : A \to A'$, compatible with the gradings. 

One defines similarly the category of group-graded co-Poisson 
bialgebras and quantization functors for group Lie bialgebras. 

\subsection{Relation with quantization of co-Poisson bialgebras}

We have inclusions of full subcategories $\{$co-Poisson universal 
enveloping algebras$\} \subset \{$co-Poisson group universal enveloping 
algebras$\} \subset \{$co-Poisson bialgebras$\}$. 

The classical limit functor is $\on{class} : \{$quasicocommutative 
bialgebras$\} \to \{$co-Poisson bialgebras$\}$. 

A quantization functor of Lie bialgebras is a functor 
$\{$co-Poisson universal enveloping algebras$\} \to \{$quasicocommutative
bialgebras$\}$, left inverse fo class.  A quantization functor for 
group Lie bialgebras may then be viewed as a left inverse to class 
with a wider domain. 

\section{The formalism of props}\label{back}

We recall material from \cite{EH}. Polynomial Schur functors form a symmetric 
monoidal abelian category Sch, equipped with an involution. 
A prop $P$ is an additive symmetric monoidal 
category, equipped with a tensor functor $\on{Sch} \to P$, which induces 
a bijection $\on{Ob(Sch)} \simeq \on{Ob}(P)$ (\cite{McL}). 
A prop morphism $P\to Q$ is a tensor functor, such that the composition 
$\on{Sch}\to P \to Q$ coincides with $\on{Sch} \to Q$. A topological 
prop is defined in the same way, with $\on{Sch}$
replaced by the category of ``formal series'' Schur functors
(i.e., infinite sums of homogenous Schur functors). 
If $F$ is a (formal series) Schur functor and $P$ is a (topological) 
prop, then $F(P)$ is 
a prop defined by $(F(P))(F_1,F_2) = P(F_1\circ F,F_2\circ F)$.  

Props may be defined by generators and relations. We will need the props 
Bialg of bialgebras, LBA of Lie bialgebras, $\on{LBA}_f$ of Lie bialgebras
with a twist. Generators of Bialg are $m,\Delta,\varepsilon,\eta$ (the 
universal analogues of the product, coproduct, counit, unit of a bialgebra); 
generators of LBA are $\mu,\delta$ (universal analogues of the Lie bracket and 
cobracket); $\on{LBA}_f$ has the additional generator $f$ (universal twist
element). LBA and $\on{LBA}_f$ are graded ($\mu$ has degree $0$ and $\delta,f$
have degree $1$) and can be completed into topological props 
${\bf LBA}$, ${\bf LBA}_f$. 

We define tensor categories $\on{Sch}_{(1)}$ and $\on{Sch}_{(1+1)}$
by $\on{Ob}(\on{Sch}_{(1)}) = \prod'_{n\geq 0} \on{Ob}(\on{Sch}_n)$
and $\on{Sch}_{(1+1)} = \prod'_{p,q\geq 0} \on{Ob}(\on{Sch}_{p+q})$, 
where $\on{Ob}(\on{Sch}_n)$ is the set of polynomial Schur 
multifunctors $\on{Vect}^n \to \on{Vect}$; the tensor product in these
categories is denoted $\boxtimes$. The bifunctor $\on{Sch}_{(1)}^2
\to \on{Sch}_{(1+1)}$ is denoted $(F,G)\mapsto F\ul\boxtimes G$. 
A multi(bi)prop
is an additive symmetric monoidal category $\tilde P$ with a tensor functor 
$\on{Sch}_{(1)} \to \tilde P$ (resp., $\on{Sch}_{(1+1)} \to \tilde P$), 
inducing a bijection on the sets of objects. A prop $P$ give rises to a 
multi-prop $\tilde P$ via $\tilde P(F,G):= P(c(F),c(G))$, where 
$c:\on{Ob}(\on{Sch}_n)\to \on{Ob}(\on{Sch})$ is induced by the diagonal 
embedding $\on{Vect}\to \on{Vect}^n$.  
We introduce the notions of a trace on a symmetric monoidal category, of a
quasi-category, we show that a symmetric monoidal category with a 
trace and an involution gives rise to a symmetric monoidal quasi-category 
(i.e., the compositions are not always defined). In particular, a trace on a 
multi-prop gives rise to a quasi-multi-bi-prop
(i.e., an additive symmetric monoidal quasi-category with a morphism from 
$\on{Sch}_{(1+1)}$ inducing a bijection on objects). 
We define traces on the multi-props arising from $\on{LBA}$ and 
$\on{LBA}_f$; this gives rise to quasi-multi-bi-props $\Pi,\Pi_f$
with $\Pi(F\underline\boxtimes G,F'\underline\boxtimes G') = 
\on{LBA}(c(F) \otimes c(G')^*,c(F')\otimes 
c(G)^*)$; the morphisms in $\Pi_f$ are defined by a similar formula. 
We also define topological completions 
\boldmath$\Pi$\unboldmath, \boldmath$\Pi$\unboldmath$_f$. 
When $F,...,G'$ are tensor products (in $\on{Sch}_{(1)}$) 
of irreducible Schur functors, 
$\Pi(F\underline\boxtimes G,F'\underline\boxtimes G')$
is graded by a set of oriented graphs; the composition of two (or several) 
morphisms is 
defined if the composition of their diagrams is acyclic. For 
general $F\in \on{Ob}(\on{Sch}_n),...,G'\in\on{Ob}(\on{Sch}_{p'})$, 
one can define the support of a given element of  
$\Pi(F\underline\boxtimes G,F'\underline\boxtimes G')$
(again an oriented graph), and acyclicity is a sufficient condition for 
the composition of two (or many) morphisms to be defined. Using this 
criterion, one checks that the compositions involved in the future 
computations all make sense. 

The motivation for working with such structures is that when $F,...,G'$
are tensor products (in $\on{Sch}_{(1)}$) of tensor Schur functors 
(i.e., objects of $\on{Sch}_1$ of the form $V\mapsto V^{\otimes n}$), 
$\Pi(F\underline\boxtimes G,F'\underline\boxtimes G')$ may be viewed as 
a space of acyclic oriented diagrams; composition is then defined by 
connecting diagrams, and is of course only defined under acyclicity
assumptions. When $F,...,G'$ are tensor products (in $\on{Sch}_{(1)}$)
of simple Schur functors, the morphisms are obtained from 
the case of tensor products of tensor functors by applying projectors 
in the group algebras of products of symmetric groups, preserving a 
partition of the vertices. 

\section{Compositions of twists} \label{comp:twists}

A quantization functor is a prop morphism $Q : \on{Bialg} \to S({\bf LBA})$ 
with certain classical limit properties. 

Let $Q$ be an Etingof-Kazhdan (EK) quantization functor. 
It is constructed as follows. 
We define elements $m_\Pi\in $\boldmath$\Pi$\unboldmath$((S\underline\boxtimes 
S)^{\boxtimes 2},
S\underline\boxtimes S)$, $\Delta_0\in $\boldmath$\Pi$\unboldmath$
(S\underline\boxtimes S,
(S\underline\boxtimes S)^{\boxtimes 2})$, 
$\on{J}\in $\boldmath$\Pi$\unboldmath$({\mathfrak 1} 
\underline\boxtimes {\mathfrak 1},
(S\underline\boxtimes S)^{\boxtimes 2})$, 
$\on{R}_+\in $\boldmath$\Pi$\unboldmath$(S\underline\boxtimes 
{\mathfrak 1},S\underline\boxtimes S)$, $m_a\in $\boldmath$\Pi$\unboldmath$
((S\underline\boxtimes {\mathfrak 1})^{\boxtimes 2}, S\underline\boxtimes
{\mathfrak 1})$, $\Delta_a\in $\boldmath$\Pi$\unboldmath$
(S\underline\boxtimes {\mathfrak 1},
(S\underline\boxtimes {\mathfrak 1})^{\boxtimes 2})$. 

We define $m_\Pi^{(i,j)}\in $\boldmath$\Pi$\unboldmath$(((S\underline
\boxtimes S
)^{\boxtimes j})^{\boxtimes i}, (S\underline\boxtimes S)^{\boxtimes j})$ 
as the $j$th tensor power of the $i$fold iterate of $m_\Pi$. 

We have 
$$
m_\Pi \circ \on{R}_+^{\boxtimes 2} = \on{R}_+ \circ m_a, \quad 
m_\Pi^{(2,2)} \circ \Big( \on{J} \boxtimes (\Delta_0 
\circ \on{R}_+)\Big) = m_\Pi^{(2,2)} \circ \Big( (\on{R}_+^{\boxtimes 2} \circ 
\Delta_a)\boxtimes \on{J}\Big).    
$$

Then $m_a,\Delta_a$ satisfy the bialgebra relations. The functor $Q$ 
is defined by $m\mapsto m_a$, $\Delta \mapsto \Delta_a$. 

We now recall the results from \cite{EH} on the quantization of twists. 
We define prop morphisms $\kappa_i : \on{LBA} \to \on{LBA}_f$ ($i = 1,2$) by 
$\kappa_1 : (\mu,\delta) \mapsto (\mu,\delta)$
and $\kappa_2 : (\mu,\delta) \mapsto (\mu,\delta + \on{ad}(f))$. 

Define $\Xi_f\in $\boldmath$\Pi$\unboldmath$_f(S\underline\boxtimes S,
S\underline\boxtimes 
S)^\times$ as the universal version of the sequence of maps 
$S(\a) \otimes S(\a^*)
\simeq U(D(\a)) \simeq U(D(\a_f)) \simeq S(\a) \otimes S(\a^*)$, 
based on the Lie algebra isomorphism $D(\a) \simeq D(\a_f)$. 
Then  
\begin{equation} \label{Xi:f}
\kappa_2^\Pi(m_\Pi) = \Xi_f \circ \kappa_1^\Pi(m_\Pi) \circ
(\Xi_f^{-1})^{\boxtimes 2}, 
\quad \kappa_2^\Pi(\Delta_0) = \Xi_f^{\boxtimes 2} \circ 
\kappa_1^\Pi(\Delta_0) \circ \Xi_f^{-1},   
\end{equation}

\begin{proposition} (see \cite{EH})
There exists $(\on{F},v,\on{i})$ with $\on{F}\in 
$\boldmath$\Pi$\unboldmath$_f({\mathfrak 1}\underline\boxtimes {\mathfrak 1},
(S\underline\boxtimes {\mathfrak 1})^{\boxtimes 2})$, $v\in 
$\boldmath$\Pi$\unboldmath$_f({\mathfrak 1}\underline\boxtimes 
{\mathfrak 1},S\underline\boxtimes S)$, and 
$\on{i}\in $\boldmath$\Pi$\unboldmath$_f(S\underline\boxtimes {\mathfrak 1}, 
S\underline\boxtimes {\mathfrak 1})^\times$, such that $\on{F} = 1
+ degree >0$, $v = 1 + degree >0$, $\on{i} = \on{id}_{S\underline\boxtimes
{\mathfrak 1}} + degree >0$, 
\begin{equation} \label{J:F}
\kappa_1^\Pi(m_\Pi^{(2,2)}) \circ 
\Big( \big( 
(\Xi_f^{-1})^{\boxtimes 2} \circ \kappa_2^\Pi(\on{J}) \big) \boxtimes 
\big( \Delta_0 \circ v\big) \Big) 
= \kappa_1^\Pi(m_\Pi^{(3,2)}) \circ \Big( \big( v\boxtimes v \big) 
\boxtimes \big( \kappa_1^\Pi(\on{R}_+)^{\boxtimes 2} \circ \on{F}\big) 
\boxtimes \big( \kappa_1^\Pi(\on{J})\big) \Big) 
\end{equation}
and 
\begin{equation} \label{id:zetaf}
\kappa_1^\Pi(m_\Pi) \circ \Big( \big( \Xi_f^{-1} 
\circ \kappa_2^\Pi(\on{R}_+) \circ \on{i} \big) \boxtimes v \Big) = 
\kappa_1^\Pi(m_\Pi) \circ \Big( v \boxtimes 
\kappa_1^\Pi(\on{R}_+) \Big). 
\end{equation}

The set of triples  $(\on{F}',v',\on{i}')$ satisfying these relations 
is given by $v' = \kappa_1^\Pi(m_\Pi) \circ \Big(v \boxtimes 
(\kappa_1^\Pi(\on{R}_+) \circ \on{u}) \Big)$, 
$m_a^{(2,2)} \circ \Big( \on{F}' \boxtimes (\Delta \circ \on{u}) \Big) 
= m_a^{(2,2)} \circ \Big( \on{u}^{\boxtimes 2} \boxtimes \on{F} \Big)$, 
$v' = \kappa_1^\Pi\Big(v \boxtimes (\kappa_1^\Pi(\on{R}_+) \circ \on{u})\Big)$, 
$\on{i}' = \on{i} \circ \underline{\on{Ad}}(\on{u})$, 
where $\on{u}\in $\boldmath$\Pi$\unboldmath$_f({\mathfrak 1} 
\underline\boxtimes {\mathfrak 1}, 
S\underline\boxtimes {\mathfrak 1})$ has the form $\on{u} = 1 + degree >0$, 
and $\underline{\on{Ad}}(\on{u})\in $\boldmath$\Pi$\unboldmath$_f
(S\underline\boxtimes {\mathfrak 1},
S\underline\boxtimes {\mathfrak 1})^\times$ is such that 
$m_a \circ (\underline{\on{Ad}}(\on{u}) \boxtimes \on{u}) = 
m_a \circ (\on{u} \boxtimes \on{id}_{S\underline\boxtimes
{\mathfrak 1}})$. 
\end{proposition}

In \cite{EH}, we prove that this proposition has the following
consequence:  

\begin{theorem} 
(Compatibility of quantization functors with twists) We have 
$$
\kappa_2^\Pi(m_a) = \on{i} \circ \kappa_1^\Pi(m_a) \circ
(\on{i}^{-1})^{\boxtimes 2}, \quad 
\kappa_1^\Pi(m_a^{(2,2)}) \circ \Big( \big( (\on{i}^{-1})^{\boxtimes 2} \circ 
\kappa_2^\Pi(\Delta_a) \circ 
\on{i}\big) \boxtimes \on{F}\Big) = 
\kappa_1^\Pi(m_a^{(2,2)}) \circ 
\Big( \on{F} \boxtimes \kappa_1^\Pi(\Delta_a)\Big),  
$$
$$
\kappa_1^\Pi(m_a^{(3,2)}) \circ \Big( (\on{F}\boxtimes 1)
\boxtimes \big( (\Delta_a \boxtimes \on{id}_{S\underline\boxtimes 
{\mathfrak 1}}) \circ \on{F}\big) \Big) = 
\kappa_1^\Pi(m_a^{(3,2)}) \circ \Big( (1\boxtimes \on{F})
\boxtimes \big( (\on{id}_{S\underline\boxtimes 
{\mathfrak 1}}\boxtimes \Delta_a) \circ \on{F}\big) \Big).  
$$
\end{theorem}

As before, $m_a^{(i,j)}$ is the $j$th tensor power of the $i$ fold 
iterate to $m_a$. 

We will now study the behavior of the composition of twists 
under quantization. 

Define a prop $\on{LBA}_{f,f'}$ by generators $\mu\in
\on{LBA}_{f,f'}(\wedge^2,{\bf id})$, $\delta\in \on{LBA}_{f,f'}({\bf id},
\wedge^2)$, $f,f'\in \on{LBA}_{f,f'}({\bf 1},\wedge^2)$ and relations: 
$\mu,\delta,f$ satisfy the relations of $\on{LBA}_f$, and $f'$ is such that 
\begin{equation} \label{rel:LBA}
((123) + (231) + (312)) \circ \Big( 
 (\delta \boxtimes \on{id}_{{\bf id}}) \circ f' +
(\mu\boxtimes \on{id}_{T_2})\circ ((1234) + (1324)) \circ (f\boxtimes f') 
\Big) =0. 
\end{equation}
Define prop morphisms  
$\kappa_{ij} : \on{LBA}_f \to \on{LBA}_{f,f'}$, by  
$\kappa_{12} : (\mu,\delta,f) \mapsto (\mu,\delta,f)$, 
$\kappa_{23} : (\mu,\delta,f) \mapsto (\mu,\delta + \on{ad}(f),f')$, 
$\kappa_{13} : (\mu,\delta,f) \mapsto (\mu,\delta,f+f')$.  

Define prop morphisms $\bar\kappa_1, \bar\kappa_2,\bar\kappa_3 : 
\on{LBA} \to \on{LBA}_{f,f'}$, $(\mu,\delta) \mapsto 
(\mu,\delta)$, $(\mu,\delta + \on{ad}(f))$, $(\mu,\delta + \on{ad}(f+f'))$. 

Then we have $\kappa_{1i} \circ \kappa_1 = \bar\kappa_1$, 
$\kappa_{i3} \circ \kappa_2 = \bar\kappa_3$,
$\kappa_{23} \circ \kappa_1 = \kappa_{12} \circ \kappa_2 = \bar\kappa_2$. 

\begin{lemma}
We have 
\begin{equation} \label{Xi:Xi}
\kappa_{23}^\Pi(\Xi_f) \circ \kappa_{12}^\Pi(\Xi_f) 
= \kappa_{13}^\Pi(\Xi_f). 
\end{equation}
\end{lemma}

{\em Proof.} This follows from the fact that if $f_\a$ is a twist for 
$\a$ and $f'_\a$ is a twist for $\a_{f_\a}$, then $f_\a + f'_\a$ is a 
twist for $\a$, and $(\a_{f_\a})_{f'_\a} \simeq (\a_{f_\a})_{f'_\a}$. 
\hfill \qed \medskip 

\begin{theorem} \label{thm:i:v} 
There exists $\on{v}\in $\boldmath$\Pi$\unboldmath$_{f,f'}({\mathfrak 1} 
\underline\boxtimes {\mathfrak 1}, S \underline{\boxtimes} {\mathfrak 1})$, 
such that $\on{v} = 1 + degree >0$,  
\begin{equation} \label{stat:v}
\bar\kappa_1^\Pi(m_a^{(2,2)}) \circ \Big( 
\kappa_{13}^\Pi(\on{F}) \boxtimes 
(\bar\kappa_1^\Pi(\Delta_a) \circ \on{v})\Big) 
= \bar\kappa_1^\Pi(m_a^{(3,2)}) \circ \Big( \on{v}^{\boxtimes 2} \boxtimes
\big( 
\kappa_{12}^\Pi(\on{i}^{-1})^{\boxtimes 2} \circ \kappa_{23}^\Pi(\on{F}) 
\big) \boxtimes \kappa_{12}^\Pi(\on{F}) \Big),  
\end{equation}
and 
\begin{equation} \label{stat:zeta}
\bar\kappa_1^\Pi(m_a) 
\circ \Big( \on{v} \boxtimes \big( \kappa_{12}^\Pi(\on{i}^{-1}) \circ 
\kappa_{23}^\Pi(\on{i}^{-1})\big) \Big)   
= \bar\kappa_1^\Pi(m_a) \circ
\Big( \kappa_{13}^\Pi(\on{i}^{-1}) \boxtimes \on{v}\Big),  
\end{equation}
\end{theorem}

{\em Proof.} Let us prove (\ref{stat:v}). 
Applying $\kappa_{23}^\Pi$ to (\ref{J:F}), we get 
\begin{align} \label{int:4}
& \bar\kappa_2^\Pi(m_\Pi^{(2,2)}) \circ 
\Big( \big( \kappa_{23}^\Pi(\Xi_f^{-1})^{\boxtimes 2} \circ 
\bar\kappa_3^\Pi(\on{J})\big) \boxtimes \big( \Delta_0 \circ 
\kappa_{23}^\Pi(v)\big) \Big) 
\\ & \nonumber
= \bar\kappa_2^\Pi(m_\Pi^{(3,2)}) \circ 
\Big( \big( \kappa_{23}^\Pi(v)^{\boxtimes 2}\big) 
\boxtimes \big( \bar\kappa_2^\Pi(\on{R}_+)^{\boxtimes 2} \circ 
\kappa_{23}(\on{F})\big) \boxtimes \big( \bar\kappa_2^\Pi(\on{J})\big) \Big).  
\end{align}
(\ref{Xi:f}) implies the identities 
\begin{equation} \label{new:id}
\bar\kappa_2^\Pi(m_\Pi) = \kappa_{12}^\Pi(\Xi_f) \circ 
\bar\kappa_1^\Pi(m_\Pi) \circ \kappa_{12}^\Pi(\Xi_f^{-1})^{\boxtimes 2}, 
\quad 
\bar\kappa_2^\Pi(\Delta_0) = \kappa_{12}^\Pi(\Xi_f)^{\boxtimes 2} \circ 
\bar\kappa_1^\Pi(\Delta_0) \circ \kappa_{12}^\Pi(\Xi_f^{-1}).  
\end{equation}
Left composing (\ref{int:4}) with $\kappa_{12}^\Pi(\Xi_f^{-1})^{\boxtimes 2}$ 
and using these identities, we get
\begin{align*}
& \bar\kappa_1^\Pi(m_\Pi^{(2,2)}) \circ 
\Big( \big( [\kappa_{12}^\Pi(\Xi_f^{-1}) \circ 
\kappa_{23}^\Pi(\Xi_f^{-1})]^{\boxtimes 2} \circ 
\bar\kappa_3^\Pi(\on{J})\big) \boxtimes \big( \Delta_0 \circ 
\kappa^\Pi_{12}(\Xi_f^{-1}) \circ \kappa_{23}^\Pi(v)\big) \Big) 
\\ & = 
\kappa^\Pi_{12}(\Xi_f^{-1})^{\boxtimes 2} \circ 
\bar\kappa_2^\Pi(m_\Pi^{(3,2)}) \circ 
\Big( \big( \kappa_{23}^\Pi(v)^{\boxtimes 2}\big) 
\boxtimes \big( \bar\kappa_2^\Pi(\on{R}_+)^{\boxtimes 2} \circ 
\kappa_{23}^\Pi(\on{F})\big) \boxtimes \big( 
\bar\kappa_2^\Pi(\on{J})\big) \Big).  
\end{align*}

Applying $\kappa_{23}^\Pi$ to (\ref{id:zetaf}), we get 
\begin{equation} \label{int:2}
\bar\kappa_2^\Pi(m_\Pi) \circ \Big( \kappa_{23}^\Pi(v) \boxtimes 
\big( \bar\kappa_2^\Pi(\on{R}_+) \big) \Big) 
= 
\bar\kappa_2^\Pi(m_\Pi) \circ \Big( 
\big( \kappa_{23}^\Pi(\Xi_f^{-1}) 
\circ \bar\kappa_3^\Pi(\on{R}_+) \circ \kappa_{23}^\Pi(\on{i})
\big) \boxtimes \kappa_{23}^\Pi(v) \Big),  
\end{equation}
and applying $\kappa_{12}^\Pi$ to (\ref{id:zetaf}), we get 
\begin{equation} \label{int:3}
\bar\kappa_1^\Pi(m_\Pi) \circ \Big( \kappa_{12}^\Pi(v) \boxtimes 
\big( \bar\kappa_1^\Pi(\on{R}_+) \big) \Big) 
= \bar\kappa_1^\Pi(m_\Pi) \circ \Big( 
\big( \kappa_{12}^\Pi(\Xi_f^{-1}) 
\circ \bar\kappa_2^\Pi(\on{R}_+) \circ \kappa_{12}^\Pi(\on{i})
\big) \boxtimes \kappa_{12}^\Pi(v) \Big).   
\end{equation}

(\ref{int:2}) then implies 
\begin{align*}
& \bar\kappa_1^\Pi(m_\Pi^{(2,2)}) \circ 
\Big( \big( [\kappa_{12}^\Pi(\Xi_f^{-1}) \circ \kappa_{23}^\Pi
(\Xi_f^{-1})]^{\boxtimes 2} \circ 
\bar\kappa_3^\Pi(\on{J})\big) \boxtimes \big( \Delta_0 \circ 
\kappa^\Pi_{12}(\Xi_f^{-1}) \circ \kappa_{23}^\Pi(v)\big) \Big) 
\\ & = 
\kappa^\Pi_{12}(\Xi_f^{-1})^{\boxtimes 2} \circ 
\bar\kappa_2^\Pi(m_\Pi^{(3,2)}) \circ 
\Big( 
\big( \kappa_{23}^\Pi(\Xi_f^{-1})^{\boxtimes 2} \circ 
\bar\kappa_3^\Pi(\on{R}_+)^{\boxtimes 2} \circ
\kappa_{23}^\Pi(\on{i})^{\boxtimes 2} \circ \kappa_{23}^\Pi(\on{F})\big) 
\boxtimes \big( \kappa_{23}^\Pi(v)^{\boxtimes 2}\big) \boxtimes 
\bar\kappa_2^\Pi(\on{J}) \Big)
\\ & = 
\bar\kappa_1^\Pi(m_\Pi^{(3,2)}) \Big( \Big( 
\big( \kappa_{12}^\Pi(\Xi_f^{-1}) \circ \kappa_{23}^\Pi(\Xi_f^{-1}) 
\circ \bar\kappa_3^\Pi(\on{R}_+) \circ \kappa_{23}^\Pi(\on{i}) 
\big)^{\boxtimes 2} \circ \kappa_{23}^\Pi(\on{F})\Big) 
\\ & 
\boxtimes \big(\kappa_{12}^\Pi(\Xi_f^{-1}) \circ \kappa_{23}^\Pi(v)
\big)^{\boxtimes 2} \boxtimes 
\big( \kappa_{12}^\Pi(\Xi_f^{-1})^{\boxtimes 2}
\circ \bar\kappa_2^\Pi(\on{J})\big) \Big) .  
\end{align*}

Applying $\kappa_{12}^\Pi$ to (\ref{J:F}), we get 
\begin{align} \label{interm}
& \nonumber \bar\kappa_1^\Pi(m_\Pi^{(2,2)}) \circ \Big( \Big( 
\kappa_{12}^\Pi(\Xi_f^{-1})^{\boxtimes 2} \circ \bar\kappa_2^\Pi(\on{J})\Big) 
\boxtimes \kappa_{12}^\Pi(\Delta_0 \circ v)\Big) 
\\ & = \bar\kappa_1^\Pi(m_\Pi^{(3,2)}) \circ \Big( \kappa_{12}^\Pi(v)^{\boxtimes 2}
\boxtimes \big( \bar\kappa_1^\Pi(\on{R}_+)^{\boxtimes 2} \circ
\kappa_{12}^\Pi(\on{F}) \big) \boxtimes \bar\kappa_1^\Pi(\on{J})\Big).  
\end{align}

Therefore "right multiplication" (using $m_\Pi$) 
of the previous identity by 
$\kappa_{12}^\Pi(\Delta_0\circ v)$ yields  
\begin{align*}
& \bar\kappa_1^\Pi(m_\Pi^{(3,2)}) \circ 
\Big( \big( [\kappa_{12}^\Pi(\Xi_f^{-1}) \circ \kappa_{23}^\Pi
(\Xi_f^{-1})]^{\boxtimes 2} \circ 
\bar\kappa_3^\Pi(\on{J})\big) \boxtimes \big( \Delta_0 \circ 
\kappa^\Pi_{12}(\Xi_f^{-1}) \circ \kappa_{23}^\Pi(v)\big) \boxtimes 
\kappa_{12}^\Pi(\Delta_0 \circ v) \Big) 
\\ & = 
\bar\kappa_1^\Pi(m_\Pi^{(4,2)}) \circ \Big( \Big( 
\big( \kappa_{12}^\Pi(\Xi_f^{-1}) \circ \kappa_{23}^\Pi(\Xi_f^{-1}) 
\circ \bar\kappa_3^\Pi(\on{R}_+) \circ \kappa_{23}^\Pi(\on{i}) 
\big)^{\boxtimes 2} \circ \kappa_{23}^\Pi(\on{F})\Big) 
\\ & 
\boxtimes \big(\kappa_{12}^\Pi(\Xi_f^{-1}) \circ \kappa_{23}^\Pi(v)
\big)^{\boxtimes 2} \boxtimes 
\big( \kappa_{12}^\Pi(\Xi_f^{-1})^{\boxtimes 2}
\circ \bar\kappa_2^\Pi(\on{J})\big) \boxtimes 
\kappa_{12}^\Pi(\Delta_0 \circ v)  \Big) 
\\ & = 
\bar\kappa_1^\Pi(m_\Pi^{(5,2)}) \circ \Big( \Big( 
\big( \kappa_{12}^\Pi(\Xi_f^{-1}) \circ \kappa_{23}^\Pi(\Xi_f^{-1}) 
\circ \bar\kappa_3^\Pi(\on{R}_+) \circ \kappa_{23}^\Pi(\on{i}) 
\big)^{\boxtimes 2} \circ \kappa_{23}^\Pi(\on{F})\Big) 
\\ & 
\boxtimes \big(\kappa_{12}^\Pi(\Xi_f^{-1}) \circ \kappa_{23}^\Pi(v)
\big)^{\boxtimes 2} \boxtimes 
\kappa_{12}^\Pi(v)^{\boxtimes 2} 
\boxtimes 
\big( \bar\kappa_1^\Pi(\on{R}_+)^{\boxtimes 2} 
\circ \kappa_{12}^\Pi(\on{F})\big) 
\boxtimes 
\bar\kappa_1^\Pi(\on{J})
\Big),   
\end{align*}
where the last equality follows from (\ref{interm}). 

According to (\ref{new:id}), the last term is equal to 
\begin{align*}
& \kappa_{12}^\Pi(\Xi_f^{-1}) \circ 
\bar\kappa_2^\Pi(m_\Pi^{(5,2)}) \circ \Big( \Big( 
\big( \kappa_{23}^\Pi(\Xi_f^{-1}) 
\circ \bar\kappa_3^\Pi(\on{R}_+) \circ \kappa_{23}^\Pi(\on{i}) 
\big)^{\boxtimes 2} \circ \kappa_{23}^\Pi(\on{F})\Big) 
\\ & 
\boxtimes \big( \kappa_{23}^\Pi(v)
\big)^{\boxtimes 2} \boxtimes 
\big(\kappa_{12}^\Pi(\Xi_f) \circ \kappa_{12}^\Pi(v) \big)^{\boxtimes 2} 
\boxtimes 
\big( \kappa_{12}^\Pi(\Xi_f)^{\boxtimes 2} \circ 
\bar\kappa_1^\Pi(\on{R}_+)^{\boxtimes 2} 
\circ \kappa_{12}^\Pi(\on{F})\big) 
\boxtimes 
\kappa_{12}^\Pi(\Xi_f)^{\boxtimes 2} \circ \bar\kappa_1^\Pi(\on{J})
\Big),   
\end{align*}
which according to (\ref{int:2}) is equal to 
\begin{align*}
& \kappa_{12}^\Pi(\Xi_f^{-1}) \circ 
\bar\kappa_2^\Pi(m_\Pi^{(5,2)}) \circ \Big( 
\kappa_{23}^\Pi(v)^{\boxtimes 2} \boxtimes 
\big( \bar\kappa_2^\Pi(\on{R}_+)^{\boxtimes 2} 
\circ \kappa_{23}^\Pi(\on{F})\big) 
\\ & 
\boxtimes 
\big(\kappa_{12}^\Pi(\Xi_f) \circ \kappa_{12}^\Pi(v) \big)^{\boxtimes 2} 
\boxtimes 
\big( \kappa_{12}^\Pi(\Xi_f)^{\boxtimes 2} \circ 
\bar\kappa_1^\Pi(\on{R}_+)^{\boxtimes 2} 
\circ \kappa_{12}^\Pi(\on{F})\big) 
\boxtimes 
\kappa_{12}^\Pi(\Xi_f)^{\boxtimes 2} \circ \bar\kappa_1^\Pi(\on{J})
\Big),   
\end{align*}
which we rewrite as 
\begin{align*}
& \bar\kappa_1^\Pi(m_\Pi^{(5,2)}) \circ \Big( 
\big( \kappa_{12}^\Pi(\Xi_f^{-1})^{\boxtimes 2} 
\circ \kappa_{23}^\Pi(v)^{\boxtimes 2}\big)  
\boxtimes 
\big( \kappa_{12}^\Pi(\Xi_f^{-1})^{\boxtimes 2} 
\circ \bar\kappa_2^\Pi(\on{R}_+)^{\boxtimes 2} 
\circ \kappa_{23}^\Pi(\on{F})\big) 
\\ & 
\boxtimes 
\big( \kappa_{12}^\Pi(v) \big)^{\boxtimes 2} 
\boxtimes 
\big( \bar\kappa_1^\Pi(\on{R}_+)^{\boxtimes 2} 
\circ \kappa_{12}^\Pi(\on{F})\big) 
\boxtimes 
\bar\kappa_1^\Pi(\on{J})
\Big).    
\end{align*}
(\ref{int:3}) allows then to rewrite this as 
\begin{align*}
& \bar\kappa_1^\Pi(m_\Pi^{(5,2)}) \circ \Big( 
\big( \kappa_{12}^\Pi(\Xi_f^{-1})^{\boxtimes 2} 
\circ \kappa_{23}^\Pi(v)^{\boxtimes 2}\big)  
\boxtimes 
\kappa_{12}^\Pi(v)^{\boxtimes 2} 
\\ & 
\boxtimes 
\big( \bar\kappa_1^\Pi(\on{R}_+)^{\boxtimes 2} \circ
\kappa_{12}^\Pi(\on{i}^{-1})^{\boxtimes 2} \circ 
\kappa_{23}^\Pi(\on{F})\big) 
\boxtimes 
\big( \bar\kappa_1^\Pi(\on{R}_+)^{\boxtimes 2} 
\circ \kappa_{12}^\Pi(\on{F})\big) 
\boxtimes 
\bar\kappa_1^\Pi(\on{J})
\Big).   
\end{align*}

We therefore get: 
\begin{equation} \label{part:1}
\bar\kappa_1^\Pi(m_\Pi^{(2,2)}) \circ \Big( \big( 
(\kappa_{12}^\Pi(\Xi_f^{-1}) \circ \kappa_{23}^\Pi(\Xi_f^{-1}))^{\boxtimes 2} 
\circ  \bar\kappa_3^\Pi(\on{J})\big) 
\boxtimes (\Delta_0 \circ v_1)\Big) = 
\bar\kappa_1^\Pi(m_\Pi^{(3,2)}) \circ
\Big( v_1^{\boxtimes 2} \boxtimes F_1 \boxtimes \bar\kappa_1^\Pi(\on{J})\Big) , 
\end{equation} 
where
$$
v_1 = \bar\kappa_1^\Pi(m_\Pi) \circ \Big( \big(
\kappa_{12}^\Pi(\Xi_f^{-1}) \circ 
\kappa_{23}^\Pi(v)\big) \boxtimes 
\kappa_{12}^\Pi(v) \Big) 
$$
\begin{align*}
F_1 & = \bar\kappa_1^\Pi(m_\Pi^{(2,2)}) \circ \Big( \big( 
\bar\kappa_1^\Pi(\on{R}_+)^{\boxtimes 2} \circ
\kappa_{12}^\Pi(\on{i}^{-1})^{\boxtimes 2} \circ \kappa_{23}^\Pi(\on{F}) 
\big) \boxtimes \big( \bar\kappa_1^\Pi(\on{R}_+)^{\boxtimes 2} \circ
\kappa_{12}^\Pi(\on{F})\big) \Big) 
\\ & = \bar\kappa_1^\Pi(\on{R}_+)^{\boxtimes 2} \circ 
\bar\kappa_1^\Pi(m_a^{(2,2)}) \circ \Big( 
\big( 
\kappa_{12}^\Pi(\on{i}^{-1})^{\boxtimes 2} \circ \kappa_{23}^\Pi(\on{F}) 
\big) \boxtimes \kappa_{12}^\Pi(\on{F}) \Big)
= \bar\kappa_1^\Pi(\on{R}_+)^{\boxtimes 2} \circ \on{F_1}, 
\end{align*}
where 
$$
\on{F}_1 = 
\bar\kappa_1^\Pi(m_a^{(2,2)}) \circ \Big( 
\big( 
\kappa_{12}^\Pi(\on{i}^{-1})^{\boxtimes 2} \circ \kappa_{23}^\Pi(\on{F}) 
\big) \boxtimes \kappa_{12}^\Pi(\on{F})
\Big). 
$$
(\ref{Xi:Xi}) implies that (\ref{part:1}) is rewritten as 
\begin{equation} \label{neweq}
\bar\kappa_1^\Pi(m_\Pi^{(2,2)}) \circ \Big( \big( 
\kappa_{13}^\Pi(\Xi_f^{-1})^{\boxtimes 2} \circ \bar\kappa_3^\Pi(\on{J})
\big) \boxtimes (\Delta_0 \circ v_1)\Big) 
= \bar\kappa_1^\Pi(m_\Pi^{(3,2)}) \circ \Big( v_1^{\boxtimes 2} 
\boxtimes \big(\bar\kappa_1^\Pi(\on{R}_+)^{\boxtimes 2} \circ \on{F}_1 \big) 
\boxtimes \bar\kappa_1^\Pi(\on{J}) \Big).
\end{equation}

On the other hand, applying $\kappa_{13}^\Pi$ to (\ref{J:F}), we get 
\begin{equation} \label{part:2}
\bar\kappa_1^\Pi(m_\Pi^{(2,2)}) \circ \Big( \big( 
\kappa_{13}^\Pi(\Xi_f^{-1})^{\boxtimes 2} 
\circ  \bar\kappa_3^\Pi(\on{J})\big) 
\boxtimes (\Delta_0 \circ v'_1)\Big) = 
\bar\kappa_1^\Pi(m_\Pi^{(3,2)}) \circ \Big( (v'_1)^{\boxtimes 2} 
\boxtimes 
\big(\bar\kappa_1^\Pi(\on{R}_+)^{\boxtimes 2} \circ \on{F}'_1 \big) 
\boxtimes \bar\kappa_1^\Pi(\on{J})\Big) , 
\end{equation} 
where 
$$
v'_1 = \kappa_{13}^\Pi(v), \quad 
\on{F}'_1 = \kappa_{13}^\Pi(\on{F}), 
$$
where $\on{F}'_1 = \kappa_{13}^\Pi(\on{F})$. 

The result of uniqueness (up to gauge) for solutions $(\on{F},v)
\in $\boldmath$\Pi$\unboldmath$_f({\mathfrak 1}\underline\boxtimes 
{\mathfrak 1},S\underline\boxtimes S)^\times \times 
$\boldmath$\Pi$\unboldmath$_f({\mathfrak 1}\underline\boxtimes 
{\mathfrak 1},S\underline\boxtimes {\mathfrak 1})^\times$
of equation (\ref{J:F}), which was established in 
\cite{EH}, Lemma 5.3, can be generalized as follows.

\begin{lemma}
The set of pairs $(\on{F}''_1,v''_1)$ satisfying (\ref{neweq}), where
$\on{F}''_1\in $\boldmath$\Pi$\unboldmath$_{f,f'}({\mathfrak 1}
\underline\boxtimes 
{\mathfrak 1},S\underline\boxtimes {\mathfrak 1})$, 
$v''_1\in $\boldmath$\Pi$\unboldmath$_{f,f'}({\mathfrak 1}\underline\boxtimes 
{\mathfrak 1},S\underline\boxtimes S)$, $\on{F}''_1 = 1 + degree >1$, 
$v''_1 = 1 + degree >1$, is given by 
$v_1 = \bar\kappa_1^\Pi(m_\Pi) \circ \Big( v''_1 \boxtimes 
\big( \bar\kappa_1^\Pi(\on{R}_+) \circ \on{v}''\big) \Big), \quad 
\bar\kappa_1^\Pi(m_a^{(2,2)}) \circ \Big( \on{F}''_1 \boxtimes 
(\bar\kappa_1^\Pi(\Delta_a) \circ \on{v}'')\Big) 
= \bar\kappa_1^\Pi(m_a^{(2,2)}) \circ \Big( \on{v}^{\prime\prime\boxtimes 2} 
\boxtimes
\on{F}_1\Big)$, where $\on{v}'' \in \Pi_{f,f'}({\mathfrak
1}\underline\boxtimes {\mathfrak 1}, S \underline\boxtimes {\mathfrak 1})$, 
$\on{v}'' = 1 +$ degree $>0$.  
\end{lemma}

The proof if parallel to that of \cite{EH}, Lemma 5.3. The computation of the 
co-Hochschild cohomology of $({\bf U}_{n,f})_{n\geq 0}$ is replaced by that of
$({\bf U}_{n,f,f'})_{n\geq 0}$, where ${\bf U}_{f,f',n} = 
\Pi_{f,f'}({\mathfrak 1}\underline\boxtimes 
{\mathfrak 1},(S\underline\boxtimes S)^{\boxtimes n})$, 
and the argument of the vanishing of $\on{LBA}_{f}({\bf id},{\bf 1})$
is replaced by the vanishing of $\on{LBA}_{f,f'}({\bf id},{\bf 1})$.  

It follows that there exists $\on{v}\in $\boldmath$\Pi$\unboldmath$_{f,f'}
({\mathfrak
1}\underline\boxtimes {\mathfrak 1}, S \underline\boxtimes {\mathfrak 1})$, 
$\on{v} = 1 +$ degree $>0$, such that 
\begin{equation} \label{final:id}
v_1 = \bar\kappa_1^\Pi(m_\Pi) \circ \Big( v'_1 \boxtimes 
\big( \bar\kappa_1^\Pi(\on{R}_+) \circ \on{v}\big) \Big), \quad 
\bar\kappa_1^\Pi(m_a^{(2,2)}) \circ \Big( \on{F}'_1 \boxtimes 
(\bar\kappa_1^\Pi(\Delta_a) \circ \on{v})\Big) 
= \bar\kappa_1^\Pi(m_a^{(2,2)}) \circ \Big( \on{v}^{\boxtimes 2} \boxtimes
\on{F}_1\Big). 
\end{equation}
The second of these identities is (\ref{stat:v}).

Let us now prove (\ref{stat:zeta}). Right composing (\ref{id:zetaf})
with $\on{i}^{-1}$, applying $\kappa_{23}^\Pi$, left composing with 
$\kappa_{12}^\Pi(\Xi_f^{-1})$, and right multiplying the resulting identity 
by $\kappa_{12}^\Pi(v)$ using $\bar\kappa_1^\Pi(m_\Pi)$, we get 
\begin{align} \label{al:1}
& \nonumber \bar\kappa_1^\Pi(m_\Pi^{(3,1)}) \circ \Big( 
\big( \kappa_{13}^\Pi(\Xi_f^{-1}) \circ \bar\kappa_3^\Pi(\on{R}_+)\big)
\boxtimes 
\big( \kappa_{12}^\Pi(\Xi_f^{-1}) \circ \kappa_{23}^\Pi(v)\big)
\boxtimes \kappa_{12}^\Pi(v)\Big)
\\ & = \bar\kappa_1^\Pi(m_\Pi^{(3,1)}) \circ \Big( 
\big( \kappa_{12}^\Pi(\Xi_f^{-1}) \circ \kappa_{23}^\Pi(v)\big)
\boxtimes \big( \kappa_{12}^\Pi(\Xi_f^{-1}) \circ \bar\kappa_2^\Pi(\on{R}_+)
\circ \kappa_{23}^\Pi(\on{i}^{-1})\big)
\boxtimes \kappa_{12}^\Pi(v)\Big). 
\end{align}

Right composing (\ref{id:zetaf}) by $\on{i}^{-1}$, applying $\kappa_{12}^\Pi$, 
right composing with $\kappa_{23}^\Pi(\on{i}^{-1})$, and left multiplying by  
$\kappa_{12}^\Pi(\Xi_f^{-1}) \circ \kappa_{23}^\Pi(v)$ using 
$\bar\kappa_1^\Pi(m_\Pi)$, we get 
\begin{align} \label{al:2}
& \nonumber \bar\kappa_1^\Pi(m_\Pi^{(3,1)}) \circ \Big( 
\big( \kappa_{12}^\Pi(\Xi_f^{-1}) \circ \kappa_{23}^\Pi(v)\big)
\boxtimes \big( \kappa_{12}^\Pi(\Xi_f^{-1}) \circ \bar\kappa_2^\Pi(\on{R}_+)
\circ \kappa_{23}^\Pi(\on{i}^{-1})\big)
\boxtimes \kappa_{12}^\Pi(v)\Big)
\\ & = 
\bar\kappa_1^\Pi(m_\Pi^{(3,1)}) \circ \Big( 
\big( \kappa_{12}^\Pi(\Xi_f^{-1}) \circ \kappa_{23}^\Pi(v)\big)
\boxtimes \kappa_{12}^\Pi(v) \boxtimes 
\big( \bar\kappa_1^\Pi(\on{R}_+) \circ \kappa_{12}^\Pi(\on{i}^{-1}) \circ 
\kappa_{23}^\Pi(\on{i}^{-1})\big)\Big). 
\end{align}

The first identity of (\ref{final:id})
is rewritten as 
\begin{equation} \label{vv:vv}
\bar\kappa_1^\Pi(m_\Pi) \circ \Big( \big( \kappa_{12}^\Pi(\Xi_f^{-1}) \circ 
\kappa_{23}^\Pi(v) \big) \boxtimes \kappa_{12}^\Pi(v) \Big)
= \bar\kappa_1^\Pi(m_\Pi) \circ \Big( \kappa_{13}^\Pi(v) \boxtimes 
(\bar\kappa_1^\Pi(\on{R}_+) \circ \on{v})\Big),   
\end{equation}
therefore 
\begin{align} \label{al:3}
& \nonumber \bar\kappa_1^\Pi(m_\Pi^{(3,1)}) \circ \Big( 
\big( \kappa_{12}^\Pi(\Xi_f^{-1}) \circ \kappa_{23}^\Pi(v)\big)
\boxtimes \kappa_{12}^\Pi(v) \boxtimes 
\big( \bar\kappa_1^\Pi(\on{R}_+) \circ \kappa_{12}^\Pi(\on{i}^{-1}) \circ 
\kappa_{23}^\Pi(\on{i}^{-1})\big)\Big). 
\\ & 
= \bar\kappa_1^\Pi(m_\Pi^{(3,1)}) \circ \Big(
\kappa_{13}^\Pi(v) \boxtimes \big( \bar\kappa_1^\Pi(\on{R}_+)\circ \on{v} \big) 
\boxtimes  
\big( \bar\kappa_1^\Pi(\on{R}_+) \circ \kappa_{12}^\Pi(\on{i}^{-1}) \circ 
\kappa_{23}^\Pi(\on{i}^{-1})\big)\Big). 
\end{align}

Combining (\ref{al:1}), (\ref{al:2}) and (\ref{al:3}), we get: 
\begin{align} \label{al:4}
& \nonumber 
\bar\kappa_1^\Pi(m_\Pi^{(3,1)}) \circ \Big( 
\big( \kappa_{13}^\Pi(\Xi_f^{-1}) \circ \bar\kappa_3^\Pi(\on{R}_+)\big)
\boxtimes 
\big( \kappa_{12}^\Pi(\Xi_f^{-1}) \circ \kappa_{23}^\Pi(v)\big)
\boxtimes \kappa_{12}^\Pi(v)\Big)
\\ & = 
\bar\kappa_1^\Pi(m_\Pi) \circ \Big(
\kappa_{13}^\Pi(v) \boxtimes   
\big[ \bar\kappa_1^\Pi(\on{R}_+) \circ \bar\kappa_1^\Pi(m_a) 
\circ \Big( \on{v} \boxtimes \big( \kappa_{12}^\Pi(\on{i}^{-1}) \circ 
\kappa_{23}^\Pi(\on{i}^{-1})\big) \Big) \big] \Big). 
\end{align}

On the other hand, left multiplying (\ref{vv:vv}) by 
$\kappa_{13}^\Pi(\Xi_f^{-1}) \circ \bar\kappa_3^\Pi(\on{R}_+)$
using $\bar\kappa_1^\Pi(m_\Pi)$, we get 
\begin{align} \label{all:1}
& \nonumber 
\bar\kappa_1^\Pi(m_\Pi^{(3,1)}) \circ \Big( 
\big( \kappa_{13}^\Pi(\Xi_f^{-1}) \circ \bar\kappa_3^\Pi(\on{R}_+)\big)
\boxtimes 
\big( \kappa_{12}^\Pi(\Xi_f^{-1}) \circ \kappa_{23}^\Pi(v)\big)
\boxtimes \kappa_{12}^\Pi(v)\Big)
\\ & 
= \bar\kappa_1^\Pi(m_\Pi^{(3,1)}) \circ \Big( 
\big( \kappa_{13}^\Pi(\Xi_f^{-1}) \circ \bar\kappa_3^\Pi(\on{R}_+)\big)
\boxtimes 
\kappa_{13}^\Pi(v) \boxtimes (\bar\kappa_1^\Pi(\on{R}_+) \circ \on{v}) \Big). 
\end{align}
Right composing (\ref{id:zetaf}) by $\on{i}^{-1}$, applying $\kappa_{13}^\Pi$
and right multiplying by $\bar\kappa_1^\Pi(\on{R}_+) \circ \on{v}$
using $\bar\kappa_1^\Pi(m_\Pi)$, we get 
\begin{align} \label{all:2}
& \nonumber 
\bar\kappa_1^\Pi(m_\Pi^{(3,1)}) \circ \Big( 
\big( \kappa_{13}^\Pi(\Xi_f^{-1}) \circ \bar\kappa_3^\Pi(\on{R}_+)\big)
\boxtimes 
\kappa_{13}^\Pi(v) \boxtimes (\bar\kappa_1^\Pi(\on{R}_+) \circ \on{v}) \Big)
\\ & \nonumber = 
\bar\kappa_1^\Pi(m_\Pi^{(3,1)}) \circ \Big( \kappa_{13}^\Pi(v) \boxtimes 
\big( \bar\kappa_1^\Pi(\on{R}_+) \circ \kappa_{13}^\Pi(\on{i}^{-1})\big) 
\boxtimes (\bar\kappa_1^\Pi(\on{R}_+) \circ \on{v})\Big)
\\ & 
= \bar\kappa_1^\Pi(m_\Pi) \circ \Big( \kappa_{13}^\Pi(v) \boxtimes 
\big[ \bar\kappa_1^\Pi(\on{R}_+) \circ \bar\kappa_1^\Pi(m_a) \circ
\Big( \kappa_{13}^\Pi(\on{i}^{-1}) \boxtimes \on{v}\Big) \big] \Big). 
\end{align}
Combining (\ref{al:4}), (\ref{all:1}) and (\ref{all:2}), we get 
\begin{align*}
& \bar\kappa_1^\Pi(m_\Pi) \circ \Big(
\kappa_{13}^\Pi(v) \boxtimes   
\big[ \bar\kappa_1^\Pi(\on{R}_+) \circ \bar\kappa_1^\Pi(m_a) 
\circ \Big( \on{v} \boxtimes \big( \kappa_{12}^\Pi(\on{i}^{-1}) \circ 
\kappa_{23}^\Pi(\on{i}^{-1})\big) \Big) \big] \Big) 
\\ & = \bar\kappa_1^\Pi(m_\Pi) \circ \Big( \kappa_{13}^\Pi(v) \boxtimes 
\big[ \bar\kappa_1^\Pi(\on{R}_+) \circ \bar\kappa_1^\Pi(m_a) \circ
\Big( \kappa_{13}^\Pi(\on{i}^{-1}) \boxtimes \on{v}\Big) \big] \Big).  
\end{align*}
Since $v$ is invertible for $m_\Pi$, and $\on{R}_+$ 
is left invertible, this implies 
$$
\bar\kappa_1^\Pi(m_a) 
\circ \Big( \on{v} \boxtimes \big( \kappa_{12}^\Pi(\on{i}^{-1}) \circ 
\kappa_{23}^\Pi(\on{i}^{-1})\big) \Big)   
= \bar\kappa_1^\Pi(m_a) \circ
\Big( \kappa_{13}^\Pi(\on{i}^{-1}) \boxtimes \on{v}\Big),  
$$
i.e., (\ref{stat:zeta}). 
\hfill \qed \medskip 

We define the prop $\on{LBA}_{f,f',f''}$ by generators 
$\mu,\delta,f,f',f''$, where $\mu\in\on{LBA}_{f,f',f''}(\wedge^2,{\bf id})$, 
$\delta\in\on{LBA}_{f,f',f''}({\bf id},\wedge^2)$, 
$f,f',f''\in\on{LBA}_{f,f',f''}({\bf 1},\wedge^2)$ and relations: 
$\mu,\delta,f,f'$ satisfy the relations of $\on{LBA}_{f,f'}$, and 
$(\mu,\delta,f+f',f'')$ satisfy the relation (\ref{rel:LBA}) satisfied
by $(\mu,\delta,f,f')$. 

Define prop morphisms 
$\kappa_{ijk} : \on{LBA}_{f,f'}
\to \on{LBA}_{f,f',f''}$, $\kappa_{123} : (\mu,\delta,f,f') \mapsto 
(\mu,\delta,f,f')$, $\kappa_{124} : (\mu,\delta,f,f') \mapsto 
(\mu,\delta,f,f'+f'')$, $\kappa_{134} : (\mu,\delta,f,f') \mapsto 
(\mu,\delta,f+f',f'')$, $\kappa_{234} : (\mu,\delta,f,f') \mapsto 
(\mu,\delta + \on{ad}(f),f',f'')$. 

Define prop morphisms $\bar{\kappa}_{ij} : \on{LBA}_f \to \on{LBA}_{f,f',f''}$
by $\bar{\kappa}_{12} : (\mu,\delta,f) \mapsto (\mu,\delta,f)$, 
$\bar{\kappa}_{13} : (\mu,\delta,f) \mapsto (\mu,\delta,f+f')$, 
$\bar{\kappa}_{14} : (\mu,\delta,f) \mapsto (\mu,\delta,f+f'+f'')$, 
$\bar{\kappa}_{23} : (\mu,\delta,f) \mapsto (\mu,\delta + \on{ad}(f),f')$, 
$\bar{\kappa}_{24} : (\mu,\delta,f) \mapsto (\mu,\delta+\on{ad}(f),
f'+f'')$, 
$\bar{\kappa}_{34} : (\mu,\delta,f) \mapsto 
(\mu,\delta+\on{ad}(f+f'),f'')$.  

Define prop morphisms $\bar{\bar{\kappa}}_i : \on{LBA} \to \on{LBA}_{f,f',f''}$ 
for $i = 1,...,4$, by $\bar{\bar{\kappa}}_1 : (\mu,\delta) \mapsto 
(\mu,\delta)$, $\bar{\bar{\kappa}}_2 : (\mu,\delta) \mapsto 
(\mu,\delta + \on{ad}(f))$, 
$\bar{\bar{\kappa}}_3 : (\mu,\delta) \mapsto 
(\mu,\delta + \on{ad}(f+f'))$, 
$\bar{\bar{\kappa}}_4 : (\mu,\delta) \mapsto 
(\mu,\delta + \on{ad}(f+f'+f''))$.  

Then $\bar\kappa_{ij} \circ \kappa_1 = \bar{\bar{\kappa}}_i$, 
$\bar\kappa_{ij} \circ \kappa_2 = \bar{\bar{\kappa}}_j$ (where 
$1\leq i<j\leq 4$); we also have 
$\kappa_{1ij} \circ \bar\kappa_1 = \bar{\bar{\kappa}}_1$ (where $2\leq 
i<j\leq 4$),
$\kappa_{234} \circ \bar\kappa_1 = \kappa_{12i} \circ \bar\kappa_2 = 
\bar{\bar{\kappa}}_2$ (where $i=3,4$),
$\kappa_{i34} \circ \bar\kappa_2 = \kappa_{123} \circ \bar\kappa_3 = 
\bar{\bar{\kappa}}_3$ (where $i=1,2$),
$\kappa_{ij4} \circ \bar\kappa_3 = \bar{\bar{\kappa}}_4$ (where 
$1\leq i<j\leq 3$); finally $\kappa_{12i} \circ \kappa_{12} = \bar\kappa_{12}$
($i=3,4$), $\kappa_{1i4} \circ \kappa_{13} = \bar\kappa_{14}$ ($i=2,3$), 
$\kappa_{i34} \circ \kappa_{23} = \bar\kappa_{34}$ ($i=1,2$), 
$\kappa_{134}\circ \kappa_{12} = \kappa_{123} \circ \kappa_{13} 
= \bar\kappa_{13}$, 
$\kappa_{234}\circ \kappa_{12} = \kappa_{123} \circ \kappa_{23} 
= \bar\kappa_{23}$, 
$\kappa_{234}\circ \kappa_{13} = \kappa_{124} \circ \kappa_{23} 
= \bar\kappa_{24}$.  

\begin{theorem}  \label{thm:v:v} \label{thm:ids:v}
\begin{equation} \label{vv=vv}
\bar{\bar\kappa}_1(m_a) \circ 
\big( \kappa_{134}^\Pi(\on{v}) \boxtimes \kappa_{123}^\Pi(\on{v})\big)
= \bar{\bar\kappa}_1(m_a) \circ \big( 
\kappa_{124}^\Pi(\on{v}) \boxtimes 
(\bar\kappa_{12}^\Pi(\on{i}^{-1})\circ \kappa_{234}^\Pi(\on{v}))\big). 
\end{equation}
\end{theorem}

{\em Proof.} Applying $\kappa_{134}^\Pi$ to (\ref{stat:v}), we get 
$$
\bar{\bar\kappa}_1^\Pi(m_a^{(2,2)}) \circ \big( \bar\kappa_{14}^\Pi(\on{F}) 
\boxtimes \kappa_{134}^\Pi(\Delta_a \circ \on{v})\big) = 
\bar{\bar\kappa}_1^\Pi(m_a^{(2,2)}) \circ \Big(
\kappa_{134}^\Pi(\on{v})^{\boxtimes 2} \boxtimes \big(
\bar\kappa_{13}(\on{i}^{-1})^{\boxtimes 2} \circ \bar\kappa_{34}(\on{F})\big)
\boxtimes \bar\kappa_{13}^\Pi(\on{F})\Big) .  
$$
Using the fact that $\bar\kappa^\Pi_{13}(\on{F}) = \kappa^\Pi_{123} 
\circ \kappa^\Pi_{13}(\on{F})$ and the image of (\ref{stat:v})
by $\kappa_{123}^\Pi$, we get 
\begin{align*}
& \bar{\bar{\kappa}}_1(m_a^{(3,2)}) \circ \Big( \bar\kappa_{14}^\Pi(\on{F})
\boxtimes ( \bar{\bar\kappa}_1^\Pi(\Delta_a)\circ 
\kappa_{134}^\Pi(\on{v})) \boxtimes
(\bar{\bar\kappa}_1^\Pi(\Delta_a) \circ \kappa_{123}^\Pi(\on{v}))\Big) 
\\ & = 
\bar{\bar{\kappa}}_1(m_a^{(5,2)}) \circ \Big( 
\kappa^\Pi_{134}(\on{v})^{\boxtimes 2} \boxtimes 
\big( \bar\kappa_{13}^\Pi(\on{i}^{-1})^{\boxtimes 2} \circ 
\bar\kappa_{34}^\Pi(\on{F})\big) 
\boxtimes \kappa_{123}^\Pi(\on{v}^{\boxtimes 2}) \boxtimes 
\big( \bar\kappa_{12}^\Pi(\on{i}^{-1})^{\boxtimes 2} \circ 
\bar\kappa_{23}^\Pi(\on{F})\big) \boxtimes \bar\kappa_{12}^\Pi(\on{F})
\Big).  
\end{align*}

Applying $\kappa_{123}^\Pi$ to (\ref{stat:zeta}), we get 
$\bar{\bar\kappa}_1^\Pi(m_a) \circ ( \bar\kappa_{13}^\Pi(\on{i}^{-1})
\boxtimes \kappa_{123}^\Pi(\on{v})) = \bar{\bar\kappa}_1^\Pi(m_a) 
\circ (\kappa_{123}^\Pi(\on{v}) \boxtimes (\bar\kappa_{12}^\Pi(\on{i}^{-1})
\circ \bar\kappa_{23}^\Pi(\on{i}^{-1})))$, which implies that 
\begin{align} \label{partial:1}
& \bar{\bar{\kappa}}_1(m_a^{(2,2)}) \circ \Big( \bar\kappa_{14}^\Pi(\on{F})
\boxtimes (\bar{\bar\kappa}_1^\Pi(\Delta_a) \circ \on{v}_1)\Big) 
\\ & = \nonumber 
\bar{\bar{\kappa}}_1(m_a^{(4,2)}) \circ \Big( 
\on{v}_1^{\boxtimes 2} \boxtimes 
\big( \bar\kappa_{12}^\Pi(\on{i}^{-1})^{\boxtimes 2} \circ 
\bar\kappa_{23}^\Pi(\on{i}^{-1})^{\boxtimes 2} \circ 
\bar\kappa_{34}^\Pi(\on{F})\big) 
\boxtimes 
\big( \bar\kappa_{12}^\Pi(\on{i}^{-1})^{\boxtimes 2} \circ 
\bar\kappa_{23}^\Pi(\on{F})\big) \boxtimes \bar\kappa_{12}^\Pi(\on{F})
\Big),   
\end{align}
where $\on{v}_1 = \bar{\bar\kappa}_1(m_a^{(2,1)}) \circ 
\big( \kappa_{134}^\Pi(\on{v}) \boxtimes \kappa_{123}^\Pi(\on{v})\big)$. 

Applying $\kappa_{124}^\Pi$ to (\ref{stat:v}), we get 
\begin{equation} \label{interm:5}
\bar{\bar\kappa}_1^\Pi(m_a)^{(2,2)} \circ \Big( \bar\kappa_{14}^\Pi(\on{F}) 
\boxtimes (\bar{\bar\kappa}_1^\Pi(\Delta_a) \circ \kappa_{124}^\Pi(\on{v}))
\Big) = 
\bar{\bar\kappa}_1^\Pi(m_a^{(3,2)}) \circ \Big( \kappa_{124}^\Pi
(\on{v}^{\boxtimes 2}) 
\boxtimes \big( \bar\kappa_{12}(\on{i}^{-1})^{\boxtimes 2} \circ
\bar\kappa_{24}(\on{F})\big) \boxtimes \bar\kappa_{12}(\on{F})\Big),  
\end{equation} 
and applying $\kappa_{234}^\Pi$ to the same identity, we get 
$$
\bar{\bar\kappa}_2^\Pi(m_a^{(2,2)}) \circ \Big( \bar\kappa_{24}(\on{F}) 
\boxtimes
(\bar{\bar\kappa}_2^\Pi(\Delta_a) \circ \kappa_{234}^\Pi(\on{v}))\Big)
= \bar{\bar\kappa}_2^\Pi(m_a^{(3,2)}) \circ \Big( \kappa_{234}
(\on{v}^{\boxtimes 2})
\boxtimes \big( \bar\kappa_{23}(\on{i}^{-1})^{\boxtimes 2} \circ 
\bar\kappa_{34}(\on{F})\big) \boxtimes \bar\kappa_{23}(\on{F})\Big).   
$$
Since $\bar{\bar\kappa}_2(m_a) = \bar\kappa_{12}(\on{i}) \circ
\bar{\bar\kappa}_1(m_a) \circ (\bar\kappa_{12}(\on{i})^{\boxtimes 2})^{-1}$, 
we get 
\begin{align*}
& \bar{\bar\kappa}_1^\Pi(m_a^{(2,2)}) \circ \Big( 
(\bar\kappa_{12}(\on{i}^{-1})^{\boxtimes 2} \circ 
\bar\kappa_{24}(\on{F})) \boxtimes
(\bar\kappa_{12}(\on{i}^{-1})^{\boxtimes 2} \circ 
\bar{\bar\kappa}_2^\Pi(\Delta_a)\circ \kappa_{234}^\Pi(\on{v}))\Big)
\\ & = \bar{\bar\kappa}_1^\Pi(m_a^{(3,2)}) \circ \Big( 
(\bar\kappa_{12}(\on{i}^{-1})^{\boxtimes 2} \circ 
\kappa_{234}^\Pi(\on{v}^{\boxtimes 2}))
\boxtimes \big( \bar\kappa_{12}(\on{i}^{-1})^{\boxtimes 2} \circ 
\bar\kappa_{23}(\on{i}^{-1})^{\boxtimes 2} \circ 
\bar\kappa_{34}(\on{F})\big) \boxtimes 
(\bar\kappa_{12}(\on{i}^{-1})^{\boxtimes 2} \circ 
\bar\kappa_{23}(\on{F})) \Big).     
\end{align*}
Right multiplying this identity by $\bar\kappa_{12}(\on{F})$ using 
$\bar{\bar\kappa}_1^\Pi(m_a)$, and 
using $\bar{\bar\kappa}_1^\Pi(m_a) \circ 
\big( [\bar\kappa_{12}^\Pi(\on{i}^{-1})^{\boxtimes 2} \circ 
\bar{\bar\kappa}_2^\Pi(\Delta_a)] \boxtimes \bar\kappa_{12}^\Pi(\on{F}) \big) 
= \bar{\bar\kappa}_1^\Pi(m_a) \circ 
\big( \bar\kappa_{12}^\Pi(\on{F}) \boxtimes 
[\bar{\bar\kappa}_1^\Pi(\Delta_a) \circ \bar\kappa_{12}^\Pi(\on{i}^{-1}) 
] \big)$, we get 
\begin{align*}
& \bar{\bar\kappa}_1^\Pi(m_a^{(3,2)}) \circ \Big( 
(\bar\kappa_{12}(\on{i}^{-1})^{\boxtimes 2} \circ 
\bar\kappa_{24}(\on{F})) \boxtimes
\bar\kappa_{12}^\Pi(\on{F}) \boxtimes 
\big( \bar{\bar\kappa}_1^\Pi(\Delta_a) \circ 
\bar\kappa_{12}^\Pi(\on{i}^{-1}) \circ \kappa_{234}^\Pi(\on{v})\big)
\Big)
\\ & = \bar{\bar\kappa}_1^\Pi(m_a^{(4,2)}) \circ \Big( 
(\bar\kappa_{12}(\on{i}^{-1}) \circ 
\kappa_{234}^\Pi(\on{v}))^{\boxtimes 2}
\boxtimes 
\\ & \big( \bar\kappa_{12}(\on{i}^{-1})^{\boxtimes 2} \circ 
\bar\kappa_{23}(\on{i}^{-1})^{\boxtimes 2} \circ 
\bar\kappa_{34}(\on{F})\big) \boxtimes 
(\bar\kappa_{12}(\on{i}^{-1})^{\boxtimes 2} \circ 
\bar\kappa_{23}(\on{F})) \boxtimes \bar\kappa_{12}^\Pi(\on{F})\Big).     
\end{align*}
Right multiplying (\ref{interm:5}) by $\bar{\bar\kappa}_1^\Pi(\Delta_a) \circ 
\bar\kappa_{12}^\Pi(\on{i}^{-1}) \circ \kappa_{234}^\Pi(\on{v})$ using 
$\bar{\bar\kappa}_1(m_a)$, we then get 
\begin{align} \label{partial:2}
& \bar{\bar{\kappa}}_1(m_a^{(2,2)}) \circ \Big( \bar\kappa_{14}^\Pi(\on{F})
\boxtimes (\bar{\bar\kappa}_1^\Pi(\Delta_a) \circ \on{v}_2)\Big) 
\\ & = \nonumber 
\bar{\bar{\kappa}}_1(m_a^{(4,2)}) \circ \Big( 
\on{v}_2^{\boxtimes 2} \boxtimes 
\big( \bar\kappa_{12}^\Pi(\on{i}^{-1})^{\boxtimes 2} \circ 
\bar\kappa_{23}^\Pi(\on{i}^{-1})^{\boxtimes 2} \circ 
\bar\kappa_{34}^\Pi(\on{F})\big) 
\boxtimes 
\big( \bar\kappa_{12}^\Pi(\on{i}^{-1})^{\boxtimes 2} \circ 
\bar\kappa_{23}^\Pi(\on{F})\big) \boxtimes \bar\kappa_{12}^\Pi(\on{F})
\Big),   
\end{align}
where $\on{v}_2 = \bar{\bar\kappa}_1(m_a^{(2,1)}) \circ \big( 
\kappa_{124}^\Pi(\on{v}) \boxtimes 
(\bar\kappa_{12}^\Pi(\on{i}^{-1})\circ \kappa_{234}^\Pi(\on{v}))\big)$. 

There exists a unique $\on{w}\in \Pi_{f,f',f''}({\mathfrak 1}\underline
\boxtimes {\mathfrak 1},
S \underline\boxtimes {\mathfrak 1})$, of the form $\on{w} = 1$ + degree $>0$, 
such that $\on{v}_1 = \bar{\bar\kappa}_1^\Pi(m_a) \circ (\on{w} \boxtimes 
\on{v}_2)$. Then (\ref{partial:1}) and (\ref{partial:2}) imply that 
$$
\bar{\bar\kappa}_1^\Pi(m_a^{(2,2)})  \circ \Big( 
\bar\kappa_{14}^\Pi(\on{F}) \boxtimes 
(\bar{\bar\kappa}_1^\Pi(\Delta_a) \circ \on{w}) \Big) 
= 
\bar{\bar\kappa}_1^\Pi(m_a^{(2,2)})  \circ \Big( 
\on{w}^{\boxtimes 2} \boxtimes \bar\kappa_{14}^\Pi(\on{F})\Big),   
$$
i.e., $\on{w}' := \bar\kappa_{14}(\on{i}) \circ \on{w}$ satisfies 
$(\on{w}')^{\boxtimes 2} = \bar{\bar\kappa}_4(\Delta_a) \circ \on{w}'$. 

Identity (\ref{vv=vv}) now follows from: 

\begin{proposition} \label{prop:gp:like}
If $x\in $\boldmath$\Pi$\unboldmath$_{f,f',f''}({\mathfrak 1}\underline
\boxtimes {\mathfrak 1},
S \underline\boxtimes {\mathfrak 1})$ 
is of the form $X = 1 +$ degree $>0$
and if $x^{\boxtimes 2} = \bar{\bar\kappa}_4(\Delta_a) \circ x$, then 
$x=1$. 
\end{proposition}

{\em Proof of Proposition.} $\on{LBA}_{f,f',f''}$ is equipped with a prop 
automorphism $\iota$, where $\iota^2 = \on{id}$, uniquely defined my 
$(\mu,\delta,f,f',f'') \mapsto (\mu,\delta + \on{ad}(f+f'+f''), -f'', 
-f',-f)$. Then $\iota \circ \bar{\bar\kappa}_4 = \bar{\bar\kappa}_1$. 
Set $y:= \iota^\Pi(x)$, then $\bar{\bar\kappa}_1^\Pi(y) 
= y^{\boxtimes 2}$. 

The prop $\on{LBA}_{f,f',f''}$ is equipped with a degree, 
such that $\on{deg}(\mu)=0$ and $\on{deg}(\delta) = \on{deg}(f) = \on{deg}(f')
= \on{deg}(f'') = 1$. We then decompose $y = 1 + y_1 + ...$ for this degree. 
Assume that we showed $y_1 = ... = y_{n-1}=0$. We then get: 
$y_n \boxtimes 1 + 1 \boxtimes y_n =$ the degree $n$ part of 
$\bar{\bar\kappa}_1^\Pi(\Delta_a) \circ y_n$, i.e., $=\Delta_0 \circ y_n$. 
According to the computation of the co-Hochschild cohomology of the complex 
$S^{\otimes 0} \to S \to S^{\otimes 2} \to...$ of Schur functors, we 
get $y_n\in \Pi_{f,f',f''}({\mathfrak 1}\underline\boxtimes {\mathfrak 1},
{\bf id}\underline\boxtimes {\mathfrak 1})
\subset \Pi_{f,f',f''}({\mathfrak 1}\underline\boxtimes {\mathfrak 1},
S \underline\boxtimes {\mathfrak 1})$.  

The degree $n+1$ part of the equation $\bar{\bar\kappa}_1^\Pi(\Delta_a) \circ y 
= y^{\boxtimes 2}$ then yields (degree $n+1$ part of
$\Delta_0\circ y_{n+1} + \bar{\bar\kappa}_1^\Pi(\Delta_a) \circ y_n) = 
y_{n+1} \boxtimes 1 + 1 \boxtimes y_{n+1}$.  Antisymmetrizing, we get 
$\delta \circ y_n = 0$. 

We then show: 

\begin{lemma}
The map $\on{LBA}_{f,f',f''}({\bf 1},{\bf id}) \to 
\on{LBA}_{f,f',f''}({\bf 1},\wedge^2)$, $y\mapsto \delta \circ y$
is injective. 
\end{lemma}

{\em Proof of Lemma.}  As in \cite{EH}, we will construct 
a retraction of this map. As in \cite{EH}, one shows that 
$\on{LBA}_{f,f',f''}(F,G)$ is the cokernel of 
$\on{LBA}(C\otimes D\otimes F,G) \to \on{LBA}(C\otimes F,G)$,
$x\mapsto x\circ ([(\on{id}_C \boxtimes p)\circ \Delta_C]\boxtimes 
\on{id}_F)$ where $C = S(\wedge^2\oplus \wedge^2 \oplus \wedge^2)$, 
$D = \wedge^3\oplus \wedge^3 \oplus \wedge^3$
$\Delta_C : C \to C^{\otimes 2}$ is induced by the coalgebra structure of
$S$, and $p \in \on{LBA}(C,D) = \oplus_{k\geq 0}\on{LBA}(S^k \circ 
(\wedge^2 \oplus \wedge^2 \oplus \wedge^2),\wedge^3\oplus \wedge^3
\oplus \wedge^3)$ has nonzero components for $k=1,2$ only; the 
$k=1$, this component specializes to $\wedge^3(\a)^{\oplus 3} \to 
\wedge^3(\a)^{\oplus 3}$, 
$$
(f_\a,f'_\a,f''_\a) \mapsto ( (\delta_\a\otimes \on{id}_\a)(f_\a)+ \on{c.p.}, 
(\delta_\a\otimes \on{id}_\a)(f'_\a)+ \on{c.p.}, 
(\delta_\a\otimes \on{id}_\a)(f''_\a)+ \on{c.p.}), 
$$
where c.p. means cyclic permutation, and for $k=2$ is specializes to 
$S^2(\wedge^3(\a)^{\oplus 3}) \to 
\wedge^3(\a)^{\oplus 3}$, 
\begin{align*}
& (f_\a,f'_\a,f''_\a)^{\otimes 2} \mapsto ( 
[f_\a^{12},f_\a^{13}] + \on{c.p.}, 
[f_\a^{12},f_\a^{\prime 13} + f_\a^{\prime 23}] 
+ [f_\a^{\prime 12},f_\a^{\prime 13}] 
+ \on{c.p.}, 
\\ & [f_\a^{12}+f_\a^{\prime 12},f_\a^{\prime\prime 13} + f_\a^{\prime\prime 23}] 
+ [f_\a^{\prime\prime 12},f_\a^{\prime\prime 13}] 
+ \on{c.p.}).  
\end{align*}
Since left and right compositions commute, we have a commutative diagram, 
$$
\begin{matrix}
\on{LBA}(C\otimes D,{\bf id}) & \stackrel{\delta\circ -}{\to} 
& \on{LBA}(C\otimes D,\wedge^2)\\ 
\downarrow & & \downarrow\\
\on{LBA}(C,{\bf id}) & \stackrel{\delta \circ -}{\to} & \on{LBA}(C,\wedge^2)
\end{matrix}
$$
whose vertical cokernel is the map $\on{LBA}_{f,f',f''}({\bf 1},{\bf id}) \to 
\on{LBA}_{f,f',f''}({\bf 1},\wedge^2)$, $y\mapsto \delta \circ y$. 

For any Schur functor $A$, we will construct a retraction 
$r_A : \on{LBA}(A,\wedge^2) \to \on{LBA}(A,{\bf id})$ of 
the map $\on{LBA}(A,{\bf id}) \to \on{LBA}(A,\wedge^2)$, 
such that the diagram 
\begin{equation} \label{comm:diag}
\begin{matrix}
\on{LBA}(C\otimes D,\wedge^2) & \stackrel{r_{C\otimes D}}{\to} 
& \on{LBA}(C\otimes D,{\bf id})\\ 
\downarrow & & \downarrow\\
\on{LBA}(C,\wedge^2) & \stackrel{r_C}{\to} & \on{LBA}(C,{\bf id})
\end{matrix}
\end{equation}
commutes. The vertical cokernel of this map is then the desired retraction. 

We have $\on{LBA}(A,{\bf id}) = \oplus_{Z\in\on{Irr(Sch)}} \on{LCA}(A,Z)
\otimes \on{LA}(Z,{\bf id})$. As in \cite{EH}, one shows that 
$\on{LCA}(Z,{\bf id} \otimes Z)$ is 1-dimensional, and one
constructs an element $\delta_Z \in \on{LCA}(Z,{\bf id} \otimes Z)$, 
such that the component $(Z',Z'') = ({\bf id},Z)$ of the 
map $\on{LA}(Z,{\bf id}) \stackrel{\delta \circ -}{\to} 
\on{LBA}(Z,\wedge^2) \subset \on{LBA}(Z,{\bf id}^{\otimes 2}) \simeq 
\oplus_{Z',Z''\in\on{Irr(Sch)}} \on{LCA}(Z,Z'\otimes Z'') \otimes 
\on{LA}(Z',{\bf id})\otimes \on{LA}(Z'',{\bf id})$ is $\lambda \mapsto 
\delta_Z \otimes \on{id}_{\bf id} \otimes \lambda$. 

It follows that the component $Z\mapsto (Z',Z'') = ({\bf id},Z)$ of the map 
$\oplus_{Z\in \on{Irr(Sch)}} \on{LCA}(A,Z) \otimes \on{LA}(Z,{\bf id}) \simeq 
\on{LBA}(A,{\bf id}) \stackrel{\delta \circ}{\to} 
\on{LBA}(A,\wedge^2) \subset \on{LBA}(A,{\bf id}^{\otimes 2})
\simeq \oplus_{Z',Z''\in\on{Irr(Sch)}} \on{LCA}(A,Z'\otimes Z'')
\otimes \on{LA}(Z',{\bf id}) \otimes \on{LA}(Z'',{\bf id})$ is 
$\kappa \otimes \lambda \mapsto (\lambda \circ \delta_Z) 
\otimes \on{id}_{{\bf id}} \otimes \kappa$.

Dually to \cite{EH}, we construct a retraction of the map $\lambda\mapsto 
\lambda \circ \delta_Z$; it gives rise to the section $r_A$. 
One then proves the commutativity of (\ref{comm:diag}) as in 
\cite{EH}. 

This ends the proof of the lemma, and therefore also of Proposition
\ref{prop:gp:like} and Theorem \ref{thm:ids:v}. 
 \hfill \qed \medskip

\bigskip

We now draw the consequences of the results of the previous Subsection
for the quantization of twists of Lie bialgebras. 

Let $\a = (\a,\mu_\a,\delta_\a)$ be a Lie bialgebra. 
Its quantization is $Q(\a) =
(S(\a)[[\hbar]],m(\a),\Delta(\a))$, where 
$m(\a) := m_a(\mu_\a,\hbar\delta_\a)$, $\Delta(\a) := 
\Delta_a(\mu_\a,\hbar\delta_\a)$. We set $\on{F}(\a,f_\a) := 
\on{F}(\mu_\a,\hbar\delta_\a,\hbar f_\a)$, $\on{i}(\a,f_\a)
:= \on{i}(\mu_\a,\hbar\delta_\a,\hbar f_\a)$. 

Then $\on{F}(\a,f_\a)\in Q(\a)^{\otimes 2}$ and $\on{i}(\a,f_\a) : 
Q(\a) \to Q(\a_{f_\a})$ are such that 
$$
m(\a_{f_\a}) = \on{i}(\a,f_\a)
\circ m(\a) \circ (\on{i}(\a,f_\a)^{\otimes 2})^{-1}, 
\quad 
\Delta(\a_{f_\a})
= \on{i}(\a,f_\a)^{\otimes 2} \circ \on{Ad}(\on{F}(\a,f_\a)) \circ 
\Delta(\a) \circ \on{i}(\a,f_\a)^{-1},
$$
$$
(\on{F}(\a,f_\a)\otimes 1) * 
(\Delta(\a)\otimes \on{id})(\on{F}(\a,f_\a)) =  
(1\otimes \on{F}(\a,f_\a)) * (\on{id} \otimes \Delta(\a))(\on{F}(\a,f_\a))
$$
(where the product $m(\a)$ is denoted $*$).  

Assume that $f'_\a$ is a twist of $\a_{f_\a}$. Then 
Theorem \ref{thm:i:v} implies that $\on{v}(\a,f_\a,f'_\a) := 
\on{v}(\mu_\a,\hbar\delta_\a,\hbar f_\a,\hbar f'_\a)$ satisfies 
$$
\on{F}(\a,f_\a + f'_\a) = \on{v}(\a,f_\a,f'_\a)^{\otimes 2} *  
(\on{i}(\a,f_\a)^{\otimes 2})^{-1}(\on{F}(\a_{f_\a},f'_\a)) *  
\on{F}(\a,f_\a) * \Delta(\a)(\on{v}(\a,f_\a,f'_\a))^{-1}, 
$$
$$
\on{i}(\a,f_\a + f'_\a) = \on{i}(\a_{f_\a},f'_\a) \circ \on{i}(\a,f_\a) 
\circ \on{Ad}(\on{v}(\a,f_\a,f'_\a)^{-1}) 
$$
(in both equalities, $m(\a)$ in understood; it is denoted $*$ in the first
equality). 

Finally, Theorem \ref{thm:v:v} implies that if $f''_\a$ is a twist of 
$\a_{f_\a + f'_\a}$, then 
$$
\on{v}(\a,f_\a+f'_\a,f''_\a) * \on{v}(\a,f_\a,f'_\a) = 
\on{v}(\a,f_\a,f'_\a+f''_\a) * \on{i}(\a,f_\a)^{-1}(\on{v}(\a_{f_\a},f'_\a,f''_\a)). 
$$

\section{Quantization of $\Gamma$-Lie bialgebras}\label{quantgamma}

\subsection{}

Assume that $(\a,\theta,f)$ is a $\Gamma$-Lie bialgebra. 
We construct its quantization as follows. Set $A = 
S(\a) \otimes \kk\Gamma[[\hbar]]$. 
We set $[x |\gamma] := x\otimes \gamma$, $[x\otimes x'|\gamma,\gamma']
:= (x\otimes \gamma) \otimes (x'\otimes \gamma')\in A^{\otimes 2}$. 

There are unique linear maps $m : A^{\otimes 2}
\to A$ and $\Delta : A\to A^{\otimes 2}$, such that 
$$
m : [x|\gamma][x'|\gamma'] \mapsto 
[x * \on{i}(\a,f_\gamma)^{-1}(\theta_\gamma(x')) *  
\on{v}(\a,f_\gamma,\wedge^2(\theta_\gamma)(f_{\gamma'}))^{-1}
|\gamma\gamma']
$$
$$
\Delta : [x|\gamma] \mapsto [\Delta(\a)(x) * 
\on{F}(\a,f_\gamma)^{-1}| \gamma,\gamma]. 
$$
The unit for $A$ is $[1|e]$, and the counit is the map 
$[x|\gamma] \mapsto \delta_{\gamma,e} \varepsilon(x)$ (recall that $*$ denotes
the product $m(\a)$ on $S(\a)[[\hbar]]$ or its tensor square). 

\begin{proposition}
This defines a bialgebra structure on $A$, 
quantizing the co-Poisson bialgebra structure induced by 
$(\a,\theta,f)$. 
\end{proposition}

{\em Proof.} This follows from the above relations on twists. 
\hfill \qed \medskip 

%Remark: $x\mapsto [x]$, $\gamma \mapsto [\gamma]$ the natural 
%maps $S(\a) \to S(\a) \otimes \kk\Gamma$, $\Gamma \to S(\a) \otimes \kk\Gamma$. 
%The product is such that 
%$$
%[\gamma][x] = [i(\mu,\delta,f_\gamma)^{-1}(\theta_\gamma(x))][\gamma], 
%$$
%$$
%[\gamma][\gamma'] = [v(f_\gamma,\wedge^2(\theta_\gamma)(f_{\gamma'}))^{-1}]
%[\gamma\gamma'], 
%$$
%$$
%[x][y] = [m(\mu,\delta)(x\otimes y)], 
%$$
%$
%\Delta([x]) = [\Delta(\mu,\delta)(x)], 
%$$
%$$
%\Delta([\gamma]) = [F(\mu,\delta,f_\gamma)^{-1}][\gamma \otimes \gamma]
%$$

\subsection{Propic version}

The quantization of $\Gamma$-Lie bialgebras has a propic version, which 
we now describe. 

Define $\on{LA}_\Gamma$ as the prop with generators 
$\mu\in \on{LA}(\wedge^2,{\bf id})$ and 
$\theta_\gamma\in \on{LA}({\bf id},{\bf id})^\times$, and relations:
Jacobi identity on $\mu$, $\Gamma \to \on{LA}_{\Gamma}
({\bf id},{\bf id})^\times$, $\gamma\mapsto 
\theta_\gamma$ is a group morphism, and $\wedge^2(\theta_\gamma) \circ \mu 
\circ \theta_\gamma^{-1} = \mu$. 

Define $\on{LBA}_\Gamma$ as the prop with generators $\mu\in 
\on{LBA}_\Gamma(\wedge^2,{\bf id})$, $\delta\in 
\on{LBA}_\Gamma({\bf id},\wedge^2)$, $\theta_\gamma\in 
\on{LBA}_\Gamma({\bf id},{\bf id})$ and $f_\gamma\in
\on{LBA}({\bf 1},\wedge^2)$, and relations: 
$(\mu,\delta)$ satisfy the relations of the prop LBA, 
$(\mu,(\theta_\gamma)_\gamma,(f_\gamma)_\gamma)$ 
satisfy the relations of $\on{LA}_\Gamma$; for each $\gamma\in\Gamma$,
$(\mu,\delta,f_\gamma)$ satisfy the defining relations of 
$\on{LBA}_f$, as well as $\wedge^2(\theta_\gamma) \circ \delta \circ 
\theta_\gamma^{-1} = \delta + \on{ad}(f_\gamma)$, 
and for each pair $\gamma,\gamma'\in \Gamma$, 
$f_{\gamma\gamma'} = f_\gamma + \wedge^2(\theta_\gamma) \circ f_{\gamma'}$. 

Define the prop $\on{Bialg}_\Gamma$ of $\Gamma$-bialgebras as follows. 
When $\Gamma$ is finite, in addition to 
the generators $m,\Delta,\eps,\eta$ of $\on{Bialg}$, it has generators 
$e_\gamma\in \on{Bialg}_\Gamma({\bf id},{\bf id})$, and the additional \
relations are 
$\sum_{\gamma\in\Gamma}e_\gamma = \on{id}_{{\bf id}}$, 
$e_\gamma \circ e_{\gamma'} = \delta_{\gamma\gamma'} e_{\gamma}$, 
$m \circ (e_{\gamma} \boxtimes e_{\gamma'}) = e_{\gamma\gamma'} \circ m$, 
$\Delta \circ e_{\gamma} = e_{\gamma}^{\boxtimes 2} \circ \Delta$, 
$e_{\gamma} \circ \eta = \delta_{\gamma e} \eta$, $\eps \circ e_{\gamma}
= \delta_{e \gamma} \eps$. 

In general, $\on{Bialg}_\Gamma$ is defined as follows. 
If $S$ is a set, define $\on{Sch}_S$ as the category of 
polynomial Schur functors $\on{Vect}^S \to V$ of the form 
$(V_s)_{s\in S} \to \oplus_{(Z_s)}
M_{(Z_s)} \otimes (\otimes_{s\in S} Z_s(V_s))$, 
where $Z_s$ are almost all ${\bf 1}$ (the unit Schur functor ${\bf 1}(V)
= {\bf k}$). Then $\on{Sch}_S$ is a symmetric tensor category. 
We define a $S$-prop as a symmetric tensor category $P$
together with a natural transformation $\on{Sch}_S \to P$, which 
is the identity on objects. A $S$-prop may be defined by generators
and relations. Then $\on{Bialg}_{(\Gamma)}$ is the $\Gamma$-prop 
defined by generators
$m_{\gamma,\gamma'}\in \on{Bialg}_{(\Gamma)}({\bf id}_\gamma\boxtimes 
{\bf id}_\gamma',{\bf id}_{\gamma\gamma'})$
$\Delta_\gamma\in \on{Bialg}_{(\Gamma)}({\bf id}_\gamma, {\bf
id}_\gamma^{\boxtimes 2})$
$\varepsilon\in \on{Bialg}_{(\Gamma)}({\bf id}_e, {\bf 1})$, 
$\eta\in \on{Bialg}_{(\Gamma)}({\bf 1},{\bf id}_e)$, and the relations 
derived from the finite case. The diagonal embedding $\on{Vect}
\to \on{Vect}^\Gamma$ gives rise to a functor $\Delta : \on{Sch}
\to \on{Sch}_S$, and we set $\on{Bialg}_\Gamma(F,G) := 
\on{Bialg}_{(\Gamma)}(\Delta(F),\Delta(G))$.  

Then any EK quantization functor gives rise to a prop morphism 
$\on{Bialg}_\Gamma \to S({\bf LBA}_\Gamma)^\Gamma$ with suitable classical 
limit properties. A group morphism $\Gamma \to \Gamma'$ gives rise to a 
commutative diagram 
$$
\begin{matrix}
\on{Bialg}_\Gamma & \to & S({\bf LBA}_\Gamma)^{\Gamma}\\ 
\downarrow  & & \downarrow \\
\on{Bialg}_{\Gamma'} & \to & S({\bf LBA}_{\Gamma'})^{\Gamma'} 
\end{matrix}
$$ 
We have therefore a quantization functor $\{$group Lie bialgebras$\}
\to \{$quasicocommutative group bialgebras$\}$ (where both sides are full 
subcategories of $\{$co-Poisson cocomutative bialgebras$\}$ and 
$\{$quasicocommutative bialgebras$\}$). 

\subsection{Quantization of quasitriangular $\Gamma$-Lie bialgebras}

We defined a quasitriangular $\Gamma$-Lie bialgebra as a triple
$(\a,r_\a,\theta_\a)$, where 
$(\a,r_\a)$ be a quasitriangular Lie bialgebra (i.e., (
$r_\a\in\a^{\otimes 2}$ satisfies the classical Yang-Baxter identity, 
and $t_\a := r_\a + r_\a^{21}$ is $\a$-invariant), and 
$\theta_\a : \Gamma\to \on{Aut}(\a,t_\a)$ be an action of $\Gamma$ by 
Lie algebra automorphisms of $\a$, preserving $t_\a$.   
It gives rise to a $\Gamma$-Lie bialgebra, with $\delta(x) = 
[r,x^1+x^2]$ and $f_\gamma := \theta_\gamma^{\otimes 2}(r_\a) - r_\a$. 

In that case a quantization can be constructed directly: 
we set $A = U(\a)\rtimes\Gamma[[\hbar]]$, the product is 
undeformed, and the coproduct is $\Delta(x) = \on{J}(\hbar r_\a)\Delta_0(x)
\on{J}(\hbar r_\a)^{-1}$ ($\Delta_0$ is the standard coproduct). 

Denote by $\on{qt}_\Gamma$ the prop of quasitriangular $\Gamma$-Lie bialgebras,
and by  ${\bf qt}_\Gamma$ its completion. We have a natural prop morphism 
$\on{LBA}_\Gamma \to \on{qt}_\Gamma$. We claim that the 
prop morphisms $\on{Bialg}_\Gamma \to S({\bf qt}_\Gamma)^\Gamma$
(the above direct construction) and the composed morphism $\on{Bialg}_\Gamma
\to S({\bf LBA}_\Gamma)^\Gamma \to S({\bf qt}_\Gamma)^\Gamma$
are equivalent (i.e., can be obtained from each other using an inner
automorphism of $S({\bf qt}_\Gamma)^\Gamma$). 

This is a consequence of the following statement on twists. 
Let $(\a,r_\a)$ be a quasitriangular Lie bialgebra and let 
$f_\a\in\wedge^2(\a)$ be a twist. According to \cite{EK}, there exists 
an invertible 
$\on{j}(\a,r_\a) : U(\a)[[\hbar]] \to S(\a)[[\hbar]]$, such that 
$m(\a) = \on{j}(\a,r_\a) \circ m_0 \circ (\on{j}(\a,r_\a)^{\otimes 2})^{-1}$, 
$\Delta(\a) = \on{j}(\a,r_\a)^{\otimes 2} \circ \on{Ad}(\on{J}(\hbar r_\a))
\circ \Delta_0 \circ \on{j}(\a,r_\a)^{-1}$. 
Then one proves that $\on{j}(\a,r_\a + f_\a) = 
\on{i}(\a,f_\a) \circ \on{j}(\a,r_\a) \circ \on{Ad}(\on{v}(\a,r_\a)^{-1})$, 
and $\on{J}(\a,r_\a + f_\a) = (\on{v}(\a,r_\a) \otimes \on{v}(\a,r_\a)) *  
(\on{j}(\a,r_\a)^{\otimes 2})^{-1}(\on{F}(\a,f_\a)) * \on{J}(\a,r_\a) *  
\Delta_0(\on{v}(\a,r_\a))^{-1}$ (here $*$ is the undeformed product on 
$U(\a)^{\otimes 2}[[\hbar]]$).  

\subsection{Open questions} 

Let $\a$ be a simple Lie algebra and let $\tilde{W}$ be its extended Weyl group.
One expects that the only possible quantization of $U(\a)\times \tilde{W}$
is the Majid-Soibelman algebra of $\a$. When $\a$ is a 
Mac-Moody Lie algebra, one expects that if $Q$ is any EK quantization  
functor, then $Q(\a,\tilde{W})$ is the Majid-Soibelman algebra of $\a$. 

Both statements are analogues of well-known results \cite{Dr,EK6}.

%\subsection{Homework in Montpellier (passerelle MCF $\to$ PR)} 
%
%Let $X$ be a formal space, $\Gamma$ be a group. 
%The $\Gamma$-graded algebras $\cO = \oplus_{\gamma\in \Gamma} \cO_\gamma$, 
%st $\cO_e = \cO_X$, correspond to: (a) group morphisms $\Gamma\to
%\on{Aut}(\cO_X)$, (b) cocycles $c : \Gamma^2 \to \cO_X^\times$. 
%
%Hopf algebra structures on such algebras correspond to: 
%(a) formal group structures on $X$, so $X = G$, 
%(b) $\Gamma\to \on{Aut}(\cO_X)$ is in fact $\Gamma\to \on{Aut}(G)$, 
%(c) a map $\Gamma \to \cO_G^{\otimes 2}$ satisfying conditions, modulo
%$\Delta(f)(f\otimes f)^{-1}$, $f\in \cO_G^\times$. 
%
%Let $\g$ be a Lie algebra, $\Gamma$ be a group. One 
%wants to classify all the "Hopf-Poisson" structures on 
%such algebras. They could be related to $\Gamma$-Lie bialgebras.
%
%One should a Hopf algebra $\cO_{A^*} \rtimes_c \Gamma$, 
%together with a "Poisson" structure. Since $\cO_{A^*} \simeq 
%\widehat S(\a)$, this could just come from the action $\Gamma \to 
%\on{Aut}(\a)$.  
%
%About the Poisson structure: ideally, after $U(\a) \rtimes \Gamma$ 
%is quantized into 
%$A = \oplus_{\gamma\in \Gamma} A_\gamma$, one
%takes the Gavarini part only for $A_\gamma$
%and then classical limit. 
%
%This "Poisson" structure should be multiplicative in a 
%suitable sense, and one could prove that 
%
%Moreover, our quantization work above should give the
%quantization of $\cO_{A^*} \rtimes \Gamma$, provided 
%one checks that the $F(f)$ and $v(f,f')$ have suitable 
%admissibility  properties.  

\end{document}